\documentclass[10pt,a4paper]{article}

\usepackage{amssymb}
\usepackage{amsmath}
\usepackage{latexsym}
\usepackage{amsfonts}
\usepackage{amscd}
\usepackage{bm}
\newtheorem{Th}{Theorem}%[section]
\newtheorem{Prop}{Proposition}%[section]
\newtheorem{Co}{Corollary}%[section]
\newtheorem{Lm}{Lemma}%[section]
\newtheorem{Lma}{Lemma}[section]

\newtheorem{Rm}{Remark}

\newcommand{\be}{\begin{equation}}
\newcommand{\ee}{\end{equation}}
\newcommand{\bes}{\begin{equation*}}
\newcommand{\ees}{\end{equation*}}
\newcommand{\R}{\mathbb{R}}
\newcommand{\N}{\mathbb{N}}

\def\Xint#1{\mathchoice
{\XXint\displaystyle\textstyle{#1}}%
{\XXint\textstyle\scriptstyle{#1}}%
{\XXint\scriptstyle\scriptscriptstyle{#1}}%
{\XXint\scriptscriptstyle\scriptscriptstyle{#1}}%
\!\int}
\def\XXint#1#2#3{{\setbox0=\hbox{$#1{#2#3}{\int}$ }
\vcenter{\hbox{$#2#3$ }}\kern-.6\wd0}}

\def\dashint{\Xint-}

\newcommand\res{\mathop{\hbox{\vrule height 7pt width .5pt depth 0pt
\vrule height .5pt width 6pt depth 0pt}}\nolimits}

\def\theequation{\thesection.\arabic{equation}}
\def\theTh{\Roman{section}.\arabic{Th}}
\def\theProp{\Roman{section}.\arabic{Prop}}
\def\theCo{\Roman{section}.\arabic{Co}}

\def\theRm{\Roman{section}.\arabic{Rm}}
\newcommand{\reset}{\setcounter{equation}{0}\setcounter{Th}{0}\setcounter{Prop}{0}\setcounter{Co}{0}\setcounter{Lma}{0}\setcounter{Rm}{0}}

\def\lf{\left}
\def\rg{\right}

\def\al{\alpha}
\def\la{\lambda}
\def\eps{\varepsilon}
\def\ds{\displaystyle}
\def\ov{\overline}

\def\p{\partial}
\def\pro{\pi_{\vec{n}}}
\def\bn{\vec{n}}

\def\bL{\vec{L}}

\def\bu{\vec{u}}

\def\bR{\vec{R}}

\def\bV{\vec{V}}

\def\bw{\vec{w}}

\def\bv{\vec{v}}

\def\bp{\vec{\Phi}}

\def\di{{D}}

\def\dio{D_{1}(x)}

\def\res{\mathop{\hbox{\vrule height 7pt width .5pt 
depth 0pt\vrule height .5pt width 6pt depth 0pt}}\nolimits}

\begin{document}

\date{ }

\title{Uniform Regularity Results\\ for Critical and Subcritical Surface Energies}
\author{Yann Bernard\footnote{School of Mathematics, Monash University,
3800 Victoria, Australia.}\:\:,\:Tristan Rivi\`ere\footnote{Department of Mathematics, ETH Zentrum,
8093 Z\"urich, Switzerland.}}

\maketitle
\noindent
$\textbf{Abstract:}$ {\it We establish regularity results for critical points to energies of immersed surfaces depending
on the first and the second fundamental form exclusively. These results hold for a large class of intrinsic elliptic Lagrangians which are sub-critical or critical.
 They are derived using uniform $\epsilon-$regularity estimates which do not degenerate as the Lagrangians approach the critical regime given by the Willmore integrand. }

%%%%%%%%%%%%%%%%%%%
%%%%%%%%%%%%%%%%%%%
%%%%%%%%%%%%%%%%%%%
%%%%%%%%%%%%%%%%%%%
%%%%%%%%%%%%%%%%%%%
%%%%%%%%%%%%%%%%%%%
%%%%%%%%%%%%%%%%%%%
%%%%%%%%%%%%%%%%%%%
%%%%%%%%%%%%%%%%%%%
%%%%%%%%%%%%%%%%%%%
%%%%%%%%%%%%%%%%%%%
%%%%%%%%%%%%%%%%%%%
%%%%%%%%%%%%%%%%%%%

\reset

\section{Introduction and Main Results}
Let $\bp:\Sigma\rightarrow\R^{3}$ denote the immersion of an oriented two-dimensional closed surface $\Sigma$ into $\R^3$. We denote by $A$ and $H$ the corresponding second fundamental form and mean curvature. In this paper, we study (weak) critical points of functionals of the types
\be\label{func}
\mathcal{W}_F\;:=\;\int_\Sigma F(H^2)\,d\text{vol}_g\qquad\text{and}\qquad \mathcal{E}_F\;:=\;\int_\Sigma F(|A|^2)\,d\text{vol}_g\:,
\ee
where $g:=\bp^*g_{\R^3}$ is the pull-back of the Euclidean metric onto $\Sigma$. Here, $F$ is a smooth function whose properties will be made more precise in due time. 
We shall denote by $\vec{n}$ the unit Gauss map given by $\vec{\Phi}$ and the orientation on $\Sigma$. With this notation, one has
\[
\forall\, X\,,\,Y\in T_x\Sigma\quad \quad A_x(X,Y)=- (d\vec{n}_xX,\vec{Y})_{g_{{\mathbb R}^3}}
\] 
where $\vec{Y}:=\vec{\Phi}_\ast Y$  and hence
\[
|A|^2=|d\vec{n}|^2_g\quad.
\]

When $F(H^2)=H^2$, the functional $\mathcal{W}_F$ is the Willmore energy, which has generated much interest for more than 2 centuries. The critical points of $\mathcal{W}_{H^2}$ are known as Willmore surfaces. The Willmore energy is invariant under the action of the non-compact M\"obius group of conformal transformations $\R^3$. As such, it does not satisfy the Palais-Smale condition and it's critical points are subject to concentration compactness phenomena. For these reasons, there is some interest of approximating the Willmore energy by more coercive energies of the form $\mathcal{W}_F$ or $\mathcal{E}_F$ where $t<<F(t)$ for $t>>1$ (for example $F(t)=(1+t)^q$ for $q>1$).

% In \cite{KLL}, it is shown that  for $F(A^2)=(1+A^2)^{p/2}$ with $p>2$, the corresponding energy ${\mathcal E}_F$ does satisfy a suitable Palais-Smale condition,  and moreover that weak critical points $\bp\in W^{2,p}$ of ${\mathcal E}_F$ are actually smooth. The same holds true for ${\mathcal  W}_F$ provided that
%\[
%\int_\Sigma|H|^2 \,d\text{vol}_g<8\pi-\delta
%\]
%for some $\delta>0$ (see again \cite{KLL}).

% In their seminal paper \cite{SU}, Sacks and Uhlenbeck construct harmonic maps by replacing the integrand $|du|^2$ with a functional involving instead $(1+|du|^2)^{p/2}$ and study the latter in the limit $p\searrow2$. An important step in the Sacks-Uhlenbeck procedure consists in obtaining $\varepsilon$-regularity estimates which are uniform in $p$. \\
In \cite{Ri3}, the second author introduced the notion of weak $W^{2,p}$ immersion for any $p\ge 2$. Let $g_0$ be any smooth reference metric on $\Sigma$, and set
\[
W^{2,p}_{imm}(\Sigma,{\R}^3):=\left\{
\begin{array}{c}
\vec{\Phi}\in W^{1,\infty}\cap W^{2,p}(\Sigma,{\mathbb R}^3)\quad; \mbox{ denoting }\quad g_{\vec{\Phi}}(X,Y):=(\vec{\Phi}_\ast X,\vec{\Phi}_\ast Y)_{g_{{\mathbb R}^3}}\\[3mm]
\quad\exists \ C_{\vec{\Phi}}>1\quad\mbox{ s. t. } C_{\vec{\Phi}}^{-1}\, g_0\le  g_{\vec{\Phi}}\le C_{\vec{\Phi}}\, g_0
\end{array}
\right\}
\]

The first part of the paper has to do with regularity properties of functionals of the form $\mathcal{W}_F$. The first result of this paper gives the  regularity of the critical points of $\mathcal{W}_F$ under some assumptions on $F$. 
\begin{Th}\label{th-regularity}
Assume that $F:[0,\infty)\rightarrow[0,\infty)$ is a $C^1$ function satisfying\footnote{$F'(s)$ is the derivative of $F(s)$ with respect to $s$.} :
\begin{itemize}
\item[(a)]\:\:\:$ C^{-1}\,t^{p}\,\leq\,F(t^2)\,\leq C\,(1+t^{p})$,
\item[(b)]\:\:\:$tF'(t^2)$ is a smooth and invertible function of $t$,
\item[(c)]\:\:\:$C^{-1}\,t^{p-1}\,\leq\,tF'(t^2)\,\leq C\, (1+t^{p-1})$,
%\item[(d)]\:\:\:$|t|^{p-1}\;\leq\;C|t|F'(t^2)$, for some positive constant $C$. 
\end{itemize}
where $C>1$ and $p\ge 1$. Then any critical point $\bp\in W^{2,p}(\Sigma,\R^3)$ of
\bes
\mathcal{W}_F(\bp)\;=\;\int_\Sigma F(H^2)\,d\text{vol}_g
\ees
is smooth. If one drops assumption (b) one has at least that $H\in L^\infty$.\hfill $\Box$
\end{Th}
A proof of this result has been established in the critical case $F(s)=s$ in \cite{Riv2} while a proof of theorem~\ref{th-regularity} can be found in \cite{KLL} in the subcritical case for $F(s)=(1+s)^q$ with $q>1$. In order to establish theorem~\ref{th-regularity} we are  adopting the parametric approach of \cite{Riv2} which is based on the local existence of isothermal coordinates (a fact which
holds ``uniformly'' for $p\ge 2$). In contrast, the proof in 
\cite{KLL} is based on the fact that, for $p>2$, any weak $W^{2,p}$ immersion is obviously locally a graph. This fails for $p=2$, and the proof in \cite{KLL} ``blows-up'' as $p\rightarrow 2$ .

Our next goal is to merge the case $p=2$ with the case $p>2$ in a single proof  by establishing estimates which remain controlled as $p\rightarrow 2$. This is realized in the following result
which is an $\eps$-regularity \underbar{independent} of $p\in [2,q]$ for any $q>2$.
\begin{Th}\label{th-epsreg}
Let $q\in(2,\infty)$. 
For any $m>0$, there exists a constant $\varepsilon_0(m)>0$ with the following property. If a conformal immersion $\bp\in W^{2,p}(D_2(0),\R^3)$ with $p\in[2,q]$ with conformal factor $\la$ is a critical point of 
\bes
\mathcal{W}_{p}\;:=\;\int(1+|H|^2)^{p/2}\,d\text{vol}_g\qquad\text{with}\qquad w_{p}\;:=\;\int_{D_2(0)}(1+|H|^2)^{p/2}\,d\text{vol}_g<\infty
\ees
satisfying
\bes
m\;:=\;\Vert\nabla\la\Vert_{L^{2,\infty}(D_2(0))}\;<\;\infty\qquad\text{and}\qquad\Vert\nabla\bn\Vert_{L^2(D_2(0))}\;<\;\varepsilon_0(m)\:,
\ees
then
\bes
\Vert H^p\, e^{2\la}\Vert_{L^\infty(D_{1}(0))}\;\leq\;C(m, w_p, q)\ \int_{D_2(0)}(1+|H|^2)^{p/2}\,d\text{vol}_g:,
\ees
for some positive constant $C(m, w_p, q)$ which is uniformly bounded as $m$, $w_p$ and $q$ are uniformly bounded.  
\hfill $\Box$
\end{Th}
%Using the previous result, standard bootstrap arguments give the following corollary.
%\begin{Co}\label{regularity}
%Let $p\ge 2$. Any critical point $\bp\in W^{2,p}(\Sigma,\R^3)$ of 
%\bes
%\mathcal{W}_p(\vec{\Phi})\;=\;\int_\Sigma (1+H^2)^{p/2}\,d\text{vol}_g
%\ees
%is smooth. \hfill $\Box$
%\end{Co}
\begin{Rm}
\label{rm-I.1} A uniform $\epsilon-$regularity for relaxations of the Willmore energy has been first established by the second author  in \cite{Ri4} for perturbations of the form
\[
W_\sigma(\vec{\Phi}):=\int_\Sigma|\vec{H}_{\vec{\Phi}}|^2\ d\text{vol}_{g_{\vec{\Phi}}}+\sigma^2\ \int_\Sigma (1+|\vec{H}_{\vec{\Phi}}|^2)^{2}\,d\text{vol}_{g_{\vec{\Phi}}}
\]
\hfill $\Box$
\end{Rm}
Combining theorem~\ref{th-epsreg} with the main weak compactness result of \cite{Ri3} gives the following corollary which can be seen as a ``Sacks-Uhlenbeck type'' concentration
compactness counterpart theorem for weak immersions which is one of the main contribution of the present work
\begin{Co}\label{Sacks-Uhlenbeck} Let $\al_k>0$ and $\al_k\rightarrow 0$ and let $\vec{\Phi}_k$ be a sequence of weak critical points of 
\[
{\mathcal W}_k(\vec{\Phi}):=\int_\Sigma (1+H^2)^{1+\al_k}\,d\text{vol}_g
\]
Assuming 
\[
\limsup_{k\rightarrow +\infty}{\mathcal W}_k(\vec{\Phi}_k)<+\infty\quad,
\]
and assuming that the sequence of Riemann surfaces associated to $(\Sigma, g_{\vec{\Phi}_k})$ is pre-compact in the moduli space of $\Sigma$, then, there exists a subsequence $k'$,
finitely many points $a_1\cdots a_Q\in \Sigma$ a sequence $\vec{p}_{k'}\in {\R}^3$ and a map $\vec{\Psi}_\infty\in C^\infty(\Sigma,{\R}^3)$ such that
\[
\vec{\Phi}_{k'}-\vec{p}_{k'}\longrightarrow \vec{\Psi}_\infty \quad\mbox{in }C^l_{loc}(\Sigma\setminus\{a_1\cdots a_Q\},{\R}^3)\quad \forall \,l\in {\N}
\]
moreover, either $\vec{\Psi}_\infty$ is a constant map or it is a smooth immersion of $\Sigma$ critical point of 
\[
{\mathcal W}(\vec{\Psi}):=\int_\Sigma (1+H^2)\,d\text{vol}_g\quad.
\]
\hfill $\Box$
\end{Co}
\begin{Rm}
\label{remark-1}
If one skips the assumption that $(\Sigma, g_{\vec{\Phi}_k})$ is pre-compact in the moduli space of $\Sigma$ a similar result holds true on each ``thick part'' of the limiting nodal surface. This is obtained combining 
theorem~\ref{th-epsreg} with theorem 0.3 of \cite{LR}.\hfill $\Box$
\end{Rm}

\medskip

Functionals of the type $\mathcal{E}_p$ are more complicated. To make it short the variations of ${\mathcal E}_p$ generate $p-$harmonic operator while the variations of ${\mathcal W}_p$ are generating standard Laplacians. When $F(|A|^2)=(1+|A|^2)^{p/2}$, it is shown in \cite{KLL} that critical points of $\int F\,d\text{vol}_g$ are smooth. In this paper, beyond regularity matters, we focus on estimates that remain stable in the limit $p\searrow2$. In particular, we prove the following $\eps$-regularity result. 
\begin{Th}\label{epsreg2}
There exist constants $\delta>0$ and for any $m>0$ there exists $\varepsilon_0(m)>0$ with the following property. If a conformal immersion $\bp\in W^{2,p}(D_1(0),\R^3)$  for $p\in[2,2+\delta]$ with conformal factor $\la$ is a critical point of
\bes
\mathcal{E}_{p}(\vec{\Phi})\;:=\;\int_{D_1(0)}(1+|d\vec{n}|_g^2)^{p/2}\,d\text{vol}_g\qquad\text{with}\qquad e_p\;:=\;\int_{D_1(0)}(1+|d\bn|_g^2)^{p/2}\,d\text{vol}_g
\ees
satisfying
\bes
m\;:=\;\Vert\nabla\la\Vert_{L^{2,\infty}(D_1(0))}\;<\;\infty\qquad\text{and}\qquad\Vert\nabla\bn\Vert_{L^2(D_1(0))}\;<\;\varepsilon_0(m)\:,
\ees
then there exists $q>2+\delta$ such that
\be
\label{vs-511}
\lf\||\nabla\vec{n}|\,e^{\la \,(2-p)/p}\rg\|_{L^q(D_{1/2}(0))}\le C(m,e_p)\ \  e_p^{1/p}\ \quad.
\ee
for some positive constant $C(m, e_p)$.\hfill $\Box$
\end{Th}
Combining theorem~\ref{epsreg2} this time with the main weak compactness result of \cite{Ri3} gives the following concentration
compactness  theorem for weak immersions which is one of the main contribution of the present work.
\begin{Co}\label{Sacks-Uhlenbeck-2} Let $\al_k>0$ and $\al_k\rightarrow 0$ and let $\vec{\Phi}_k$ be a sequence of weak critical points of 
\[
{\mathcal E}^k(\vec{\Phi}):=\int_\Sigma (1+|d\vec{n}|^2)^{1+\al_k}\,d\text{vol}_g
\]
Assuming 
\[
\limsup_{k\rightarrow +\infty}{\mathcal E}^k(\vec{\Phi}_k)<+\infty\quad,
\]
and assuming that the sequence of Riemann surfaces associated to $(\Sigma, g_{\vec{\Phi}_k})$ is pre-compact in the moduli space of $\Sigma$, then, there exists a subsequence $k'$,
finitely many points $a_1\cdots a_Q\in \Sigma$, a sequence $\vec{p}_{k'}\in {\R}^3$ and a map $\vec{\Psi}_\infty\in C^\infty(\Sigma,{\R}^3)$ such that
\[
\vec{\Phi}_{k'}-\vec{p}_{k'}\longrightarrow \vec{\Psi}_\infty \quad\mbox{in }C^1_{loc}(\Sigma\setminus\{a_1\cdots a_Q\},{\R}^3)\quad	\quad.
\]
Moreover, either $\vec{\Psi}_\infty$ is a constant map or it is a smooth immersion of $\Sigma$ critical point of 
\[
{\mathcal W}(\vec{\Psi}):=\int_\Sigma (1+H^2)\,d\text{vol}_g\quad.
\]\hfill $\Box$
\end{Co}
\begin{Rm}
\label{remark-2}
Here again, if one skips the assumption that $(\Sigma, g_{\vec{\Phi}_k})$ is pre-compact in the moduli space of $\Sigma$ a similar result holds true on each ``thick part'' of the limiting nodal surface. This is obtained combining 
theorem~\ref{epsreg2} with theorem 0.3 of \cite{LR}.\\[-3ex]

\hfill $\Box$
\end{Rm}
\begin{Rm}
\label{general-dimensions} 
Most of the estimates bellow extend to the corresponding tensors for immersions into Euclidian Spaces of arbitrary dimensions. Nevertheless, to make the proof more readable we chose to present
\end{Rm}

\noindent
\textbf{Acknowledgment:\:\:} This work has been initiated  while the first author was visiting the Forschungsinstitut f\"ur Mathematik at the ETH in Z\"urich.  He would like to
thank the institute for the excellent working conditions and stimulating environment.

\setcounter{equation}{0} 
\reset

\section{Preliminaries}
\reset
\subsection{First Variations of ${\mathcal W}_F$ and ${\mathcal E}_F$.}\label{vary}

In local coordinates $\{x_1,x_2\}$ on the unit disk, let the pull-back metric by $\bp$ be $g_{ij}:=\partial_{x_i}\bp\cdot\partial_{x_j}\bp$. As usual, we let $|g|$ denote the determinant of the matrix $(g_{ij})$. We let $\bn$ denote the outward-unit normal vector.
We will suppose that $\bp$ is conformal, namely, in local coordinates $\{x_1,x_2\}$, we have $\partial_{x_i}\bp\cdot\partial_{x_j}\bp=\text{e}^{2\la}\delta_{ij}$. In order to derive the Euler-Lagrange equation associated with the energies $\mathcal{W}_F$ and $\mathcal{E}_F$, we consider a variation of the type
\bes
\bp_t\;:=\;\bp+t\bw+\text{o}(t)\:,
\ees
where $\bw$ is a normal vector. The variation of the outward unit normal vector $\bn$ is easily found to be
\bes
\bn_t\;=\;\bn+t\text{e}^{-\la}(a_1\partial_{x_1}\bp+a_2\partial_{x_2}\bp)+\text{o}(t)\:,
\ees
for some $a_1$ and $a_2$. One easily verifies that
\be\label{varn}
\dfrac{d}{dt}\bigg|_{t=0}\bn_t\;=\;-\,\langle\bn\cdot d\bw\,,\,\bp\rangle_{g}\:.
\ee
The variation of the components of the metric $(g_t)_{ij}:=\partial_{x_i}\bp_t\cdot\partial_{x_j}\bp_t$ is also easily found to be
\be\label{varmet}
\dfrac{d}{dt}\bigg|_{t=0}(g_t)_{ij}\;=\;\partial_{x_i}\bp\cdot\partial_{x_j}\bw+\partial_{x_j}\bp\cdot\partial_{x_i}\bw\:.
\ee
From this and the fact that $(g_t)_{ki}(g_t)^{ij}=\delta^j_k$, it follows that the inverse metric coefficients vary according to
\be\label{varinvmet}
\dfrac{d}{dt}\bigg|_{t=0}g^{ij}_t\;=\;-\,\text{e}^{-4\la}\big(\partial_{x_i}\bp\cdot\partial_{x_j}\bw+\partial_{x_j}\bp\cdot\partial_{x_i}\bw\big)\:.
\ee
Note that the variation of the volume form is
\be\label{varvol}
\dfrac{d}{dt}\bigg|_{t=0}d\text{vol}_{g_t}\;=\;\langle d\bp\,;\,d\bw\rangle_g\,d\text{vol}_g\:.
\ee
The variations of the mean curvature and of the square of the second fundamental form satisfy respectively
\be\label{varmc}
\dfrac{d}{dt}\bigg|_{t=0}H^2_t\;=\;H\big(\langle d\bn\,;\,d\bw\rangle_g-d^{*_g}(\bn\cdot d\bw)\big)
\ee
and
\begin{eqnarray}\label{var2ff}
\dfrac{d}{dt}\bigg|_{t=0}|d\bn_t|_g^2&\equiv&\dfrac{d}{dt}\bigg|_{t=0}\big(g_t^{ij}\partial_{x_i}\bn_t\cdot\partial_{x_j}\bn_t  \big)\nonumber\\[1ex]
&=&-\,2\big\langle d\bp\stackrel{.}{\otimes}d\bw\,,\,d\bn\stackrel{.}{\otimes}d\bn\big\rangle_g-2\big\langle d\langle\bn\cdot d\bw\,,\,d\bp\rangle_g\,;\,d\bn\big\rangle_g\:,
\end{eqnarray}
where we have used (\ref{varinvmet}) and (\ref{varn}). 
Combining (\ref{varvol}) and (\ref{varmc}) gives
\begin{eqnarray*}
\dfrac{d}{dt}\bigg|_{t=0}\mathcal{W}_F&=&\int_\Sigma \Big[F(H^2)\langle d\bp\,;\,d\bw\rangle_g+HF'(H^2)\big(\langle d\bn\,;\,d\bw\rangle_g-d^{*_g}(\bn\cdot d\bw)\big)\Big]\,d\text{vol}_g\:.
\end{eqnarray*}
where $F'(H^2)$ is understood as $dF(H^2)/dH^2$. \\
It is a simple matter to integrate the latter by parts. Equating the resulting integral to zero and recalling that the obtained identity holds for all $\bw$, one concludes that a critical point of $\mathcal{W}_F$ must satisfy the Euler-Lagrange equation
\begin{eqnarray}\label{EL}
&&d^{*_g}\Big[Fd\bp+HF'd\bn-\bn\,d(HF') \Big]\;=\;0\:.
\end{eqnarray}

Combining now (\ref{varvol}) and (\ref{var2ff}) and using the fact that $|A|^2=|d\bn|_g^2$ gives
\begin{eqnarray}
\dfrac{d}{dt}\bigg|_{t=0}\mathcal{E}_F&=&\int_\Sigma \Big[F\langle d\bp\,;\,d\bw\rangle_g-\;2F'\big[\big\langle d\bp\stackrel{.}{\otimes}d\bw\,,\,d\bn\stackrel{.}{\otimes}d\bn\big\rangle_g+\big\langle d\langle\bn\cdot d\bw\,,\,d\bp\rangle_g\,;\,d\bn\big\rangle_g   \big]\Big]\,d\text{vol}_g\:,
\end{eqnarray}
where $F':=dF(|A|^2)/d|A|^2$. \\
It is a simple matter to integrate the latter by parts. Equating the resulting integral to  zero and recalling that the obtained identity holds for all $\bw$, one concludes that a critical point of $\mathcal{E}_F$ must satisfy the Euler-Lagrange equation
\be\label{hoopla}
d^{*_g}\Big[{F}d\bp-\,2F'(d\bn\stackrel{.}{\otimes}d\bn)\res_gd\bp+2\big(d^{*_g}F'd\bn\big)\cdot d\bp\big)\bn \Big]\;=\;0\:,
\ee
where $(d\bn\stackrel{.}{\otimes}d\bn)\res_gd\bp$ is the contraction given in local conformal coordinates by
\bes
(d\bn\stackrel{.}{\otimes}d\bn)\res_gd\bp\;=\;\text{e}^{-2\la}\sum_{i,j=1}^2\big(\partial_{x_i}\bn\cdot\partial_{x_j}\bn\big)\partial_{x_j}\bp\,dx_i\:.
\ees

\subsection{Conservation laws for critical points of $\mathcal{W}_F$}

\begin{Prop}
\label{pr-Euler}
In a conformal chart, a critical point of $\mathcal{W}_F$, respectively $E_F$, satisfies the Euler-Lagrange equation (\ref{EL}), which reads
\begin{eqnarray}\label{divform}
&&\text{div}\Big[F(H^2)\nabla\bp+HF'(H^2)\nabla\bn-\bn\nabla (HF'(H^2))\Big]\;=\;0\:.
\end{eqnarray}
respectively for $E_F$ the equation (\ref{hoopla}) which reads
\be\label{divee-0}
\text{div}\Big[F\nabla\bp-2\,\text{e}^{-2\la}F'\sum_{j=1}^{2}(\nabla\bn\cdot\partial_{x_j}\bn)\partial_{x_j}\bp+2\text{e}^{-2\la}\bn\big(\text{div}\big(F'\nabla\bn\big)\cdot\nabla\bp\big)  \Big]\;=\;0\:.
\ee

\hfill $\Box$
\end{Prop}

We shall now first concentrate on critical points of $\mathcal{W}_F$.
The divergence form (\ref{divform}) may be locally integrated to yield the existence of a potential function $\bL$ satisfying
\be\label{defL}
\nabla^\perp\bL\;=\;HF'(H^2)\nabla\bn-\bn\nabla\big(H F'(H^2)\big)+F(H^2)\nabla\bp\:,
\ee
where, in terms of the local coordinates $\{x_1,x_2\}$, we have set $\nabla^\perp:=(-\partial_{x_2},\partial_{x_1})$. \\
Note that
\bes
\nabla\bn\cdot\nabla\bp\;=\;-\bn\cdot\Delta\bp\;=\;-2\text{e}^{2\la}H\:.
\ees
Hence (\ref{defL}) yields
\bes
\nabla\bp\cdot\nabla^\perp\bL\;=\;2\,\text{e}^{2\la}\big(F-H^2F'  \big)\:.
\ees
For the sake of our future needs, let $Y$ be a solution of
\be\label{defY00}
-\,\Delta Y\;=\;2\text{e}^{2\la}\big(F-H^2F'  \big)\:.
\ee
%\begin{Rm}
%Observe that $Y$ may be to chosen to identically vanish if and only if $F=H^2$, that is if and only if $\mathcal{W}_p$ is the Willmore energy. This is not surprising. Indeed, unlike general energies of the type $\mathcal{W}_p$, the Willmore energy is invariant under dilation (whereas the generic members of the type $\mathcal{W}_p$ are only invariant under translation and rotation). In that case, (\ref{inter1}) becomes an exact divergence form equation which results from applying Noether's theorem to the Willmore energy for its invariance under dilation (cf. \cite{Ber} for details). Our goal in this paper is not only to produce a regularity proof for critical points of $\mathcal{W}_p$, but more importantly to exhibit an argument which is stable in the limit $p\searrow2$, i.e. when approaching the Willmore energy, which is well-known to be critical. The behavior of the function $Y$ thus plays a decisive role in our argument. 
%\end{Rm}
Equation (\ref{defY00}) states that
\bes
\text{div}\big(\bL\cdot\nabla^\perp\bp-\nabla Y\big)\;=\;0\:.
\ees
This identity is integrated to give the existence of a function $S$ satisfying
\be\label{defS}
\nabla S\;=\;\bL\cdot\nabla\bp+\nabla^\perp Y\:.
\ee
Next, using that $\nabla\bp\times\nabla\bn=\text{div}(\nabla\bp\times\bn)=\text{div}\nabla^\perp\bp=0$, we obtain from (\ref{defL}) that
\bes
\nabla\bp\times\nabla^\perp\bL\;=\;-\nabla^\perp\bp\cdot\nabla\big(HF'\big)\;=\;-\,\text{div}\big(HF'\nabla^\perp\bp \big)\:,
\ees
which is an exact divergence equation and may thus be integrated to give the existence of a potential function $\bV$ satisfying
\be\label{defV0}
\nabla\bV\;=\;\bL\times\nabla\bp+HF'\nabla\bp\:.
\ee
We can summarize what we have established so far in the following proposition.
\begin{Prop}
\label{pr-cons-law}
Let $\vec{\Phi}$ be a weak critical point of
\[
{\mathcal W}_F(\vec{\Phi}):=\int_\Sigma F(H^2)\ d\text{vol}_{g}\:.
\]
In local conformal coordinate we introduce\footnote{The local existence of $\vec{L}$ is given by Proposition~\ref{pr-Euler}.}
$\vec{L}$ such that
\[
\nabla^\perp\bL\;=\;HF'(H^2)\nabla\bn-\bn\nabla\big(H F'(H^2)\big)+F(H^2)\nabla\bp\:,
\]
Then the following two identities hold
\be
\label{cons-law}
\lf\{
\begin{array}{l}
\ds\nabla\bp\cdot\nabla^\perp\bL\;=\;2\,\text{e}^{2\la}\big(F(H^2)-H^2F'(H^2)  \big)\\[3mm]
\ds\nabla\bp\times\nabla^\perp\bL\;=\;-\,\text{div}\big(H\,F'(H^2)\nabla^\perp\bp \big)\:,
\end{array}
\rg.
\ee
where $e^\la:=|\p_{x_1}\vec{\Phi}|=|\p_{x_2}\vec{\Phi}|$.
\hfill $\Box$
\end{Prop}

We now derive central identities linking together the potentials $S$ and $\bV$. Note first that
\be\label{interr00}
\bn\cdot\nabla\bV\;=\;\bn\cdot(\bL\times\nabla\bp)\;=\;\bL\cdot(\nabla\bp\times\bn)\;=\;\bL\cdot\nabla^\perp\bp\;=\;\nabla^\perp S+\nabla Y\:.
\ee
The tangential and normal parts of $\nabla\bV$ are
\bes
\pi_T\nabla\bV\;=\;-(\bL\cdot\bn)\nabla^\perp\bp+HF'\nabla\bp\qquad\text{and}\qquad\pro\nabla\bV\;=\;\big(\nabla^\perp S+\nabla Y\big)\bn\:.
\ees
Hence
\begin{eqnarray}\label{interr01}
\bn\times\nabla\bV&\equiv&\bn\times\pi_T\nabla\bV\;\;=\;\;-(\bL\cdot\bn)\nabla\bp-HF'\nabla^\perp\bp\nonumber\\[1ex]
&=&-\,\pi_T\nabla^\perp\bV\;=\;-\,\nabla^\perp\bV-\big(\nabla S-\nabla^\perp Y\big)\bn\:.
\end{eqnarray}
We formally decompose \`a la Hodge the quantity
\be\label{hodgey}
\bn\nabla^\perp Y\;=\;\nabla\vec{v}+\nabla^\perp\vec{u}\:,
\ee
and set $\bR:=\bV-\vec{u}$ into (\ref{interr00}) and into (\ref{interr01}) to discover the following proposition
\begin{Prop}
\label{pr-euler-cons-law}
With the previous notations, the following equation hold
\be\label{cwsfirst0}
\bn\times\nabla\bR=-\,\nabla^\perp\bR-\bn\nabla S-\bn\times\nabla^\perp\vec{v}+\nabla\vec{v}\quad.
\ee
\hfill $\Box$
\end{Prop}
\subsection{Control of the conformal factor.}
Using F. H\'elein's method of moving Coulomb frames \cite{ParksCity}, a weak immersion $\bp\in W^{2,2}_{imm}(D_2(0),{\R}^3)$ of the unit disk $D_2(0)$ into $\R^3$ can be re-parametrized by a diffeomorphism of $D_2(0)$ to become conformal. Our functional being independent of parametrization, we will without loss of generality suppose that $\bp$ is conformal with parameter $\la$, namely:
\bes
\partial_{x_i}\bp\cdot\partial_{x_j}\bp\;=\;\text{e}^{2\la}\delta_{ij}\:.
\ees
We will henceforth use the notation $\nabla$, $\text{div}$, and $\Delta$ to denote the usual gradient, divergence, and Laplacian operators in flat local coordinates $\{x_1,x_2\}$. \\
Assume
\[
\int_{D_2(0)}|\nabla\vec{n}|^2\ dx^2\le 8\pi/3\quad \mbox{ and }\quad \|\nabla\la\|_{L^{2,\infty}(D^2(0))}=m<+\infty\quad.
\]
We can call upon Lemma 5.1.4 in \cite{Hel} to deduce the existence of an orthogonal frame $\{\vec{e}_1,\vec{e}_2\}\in W^{1,2}(D_1(0))$ satisfying $\bn=\vec{e}_1\times\vec{e}_2$ and
\be
\label{coulomb}
\Vert\nabla\vec{e}_1\Vert_{L^2(D_2(0))}+\Vert\nabla\vec{e}_2\Vert_{L^2(D_2(0))}\;\leq\;C\, \|\nabla\vec{n}\|_{L^2(D_2(0))}\:,
\ee
 As is easily verified, the conformal parameter satisfies
\be
\label{la-equ}
\Delta\la\;=\;\nabla\vec{e}_1\cdot\nabla^\perp\vec{e}_2\qquad\text{in}\:\:D_2(0)\:.
\ee
Let $\mu$ satisfy
\be
\label{mu-equ}
\left\{\begin{array}{rclcl}
\Delta\mu&=&\nabla\vec{e}_1\cdot\nabla^\perp\vec{e}_2&,&\text{in}\:\:D_2(0)\\[1ex]
\mu&=&0&,&\text{on}\:\:\partial D_2(0)\:.
\end{array}\right.
\ee
Standard Wente estimates (cf. Theorem 3.4.1 in \cite{Hel}) give
\be\label{wenteagain}
\Vert \mu\Vert_{L^\infty(D_2(0))}+\Vert \nabla\mu\Vert_{L^2(D_2(0))}\;\leq\;\Vert\nabla\vec{e}_1\Vert_{L^2(D_2(0))}\Vert\nabla\vec{e}_2\Vert_{L^2(D_2(0))}\;\leq\;  C\, \|\nabla\vec{n}\|^2_{L^2(D_2(0))}  \:.
\ee
The harmonic function $\nu:=\la-\mu$ satisfies the usual estimate
\bes
\int_{D_{3/2}(0)}|\nu-\bar{\nu}|\,dx^2\;\leq\;C\,\Vert\nabla\nu\Vert_{L^1(D_2(0))}\;\leq\;C\,\Vert\nabla\nu\Vert_{L^{2,\infty}(D_2(0))}\:,
\ees
where $\bar{\nu}$ denotes the average of $\nu$ on $ D_1(0)$. Hence
\bes
\Vert\nu-\bar{\nu}\Vert_{L^\infty(D_{3/2}(0))}\;\leq\;C\,\Vert\nabla\nu\Vert_{L^{2,\infty}(D_2(0))}\:,
\ees 
Combining the latter to (\ref{wenteagain}) yields now
\bes
\Vert\la-\bar{\la}\Vert_{L^\infty(D_{3/2}(x))}\;\leq\;C\,\Vert\nabla\la\Vert_{L^{2,\infty}(D_2(0))}+\;  C\, \|\nabla\vec{n}\|^2_{L^2(D_2(0))} \;\leq\;C(m)\:,
\ees
where $\bar{\la}$ denotes the average of $\la$ on $D_{3/2}(0)$, and $m$ is the $L^{2,\infty}$ norm of $\nabla\la$ on $D_2(0)$, assumed to be finite. We can summarize this subsection
by stating the following lemma.
\begin{Lm}
\label{conf-factor}
Let $\bp\in W^{2,2}_{imm}(D_2(0),{\R}^3)$ be a conformal weak immersion such that
\[
\int_{D_2(0)}|\nabla\vec{n}|^2\ dx^2\le 8\pi/3\quad \mbox{ and }\quad \|\nabla\la\|_{L^{2,\infty}(D^2(0))}=m<+\infty\quad.
\]
where $e^\la:=|\p_{x_1}\vec{\Phi}|=|\p_{x_2}\vec{\Phi}|$. Then the following estimate holds
\be
\label{conf-factor-control}
\|\nabla\la\|_{L^2(D_{3/2}(0))}+\lf|\la-\frac{1}{|D_{3/2}(0)|}\int_{D_{3/2}(0)}\la(x)\ dx^2\rg|\le \;C\,\Vert\nabla\la\Vert_{L^{2,\infty}(D_2(0))}+\;  C\, \|\nabla\vec{n}\|^2_{L^2(D_2(0))}\le\, C(m) \quad.
\ee
\hfill $\Box$
\end{Lm}
We shall now prove the following extension to general exponents $p\in [2,+\infty)$.
\begin{Lm}
\label{lm-conf-factor-Lp}
Let $2\le p<q$ $\bp\in W^{2,p}_{imm}(D_2(0),{\R}^3)$ be a conformal weak immersion such that
\[
\int_{D_2(0)}|\nabla\vec{n}|^2\ dx^2\le 8\pi/3\quad \mbox{ and }\quad \|\nabla\la\|_{L^{2,\infty}(D^2(0))}=m<+\infty\quad.
\]
where $e^\la:=|\p_{x_1}\vec{\Phi}|=|\p_{x_2}\vec{\Phi}|$. Then the following estimate holds
\be
\label{conf-factor-control-Lp}
\|\nabla\la\|_{L^p(D_{3/2}(0))}\le \;C(q)\,\Vert\nabla\la\Vert_{L^{2,\infty}(D_2(0))}+\;  C(q)\, \|\nabla\vec{n}\|_{L^p(D_2(0))} \quad.
\ee
where $C(q)>0$ is independent of $p\in[2,q]$.
\hfill $\Box$
\end{Lm}
\noindent{\bf Proof of lemma~\ref{lm-conf-factor-Lp}.} The first step consists in constructing an orthonormal frame $(\vec{e}_1,\vec{e}_2)$ in $W^{1,p}(D_2(0))$ satisfying
$\vec{n}=\vec{e}_1\times\vec{e}_2$ and 
\be
\label{conf-fact-lp-1}
\begin{array}{l}
\ds \Vert\nabla\vec{e}_1\Vert_{L^2(D_2(0))}+\Vert\nabla\vec{e}_2\Vert_{L^2(D_2(0))}\;\leq\;C\, \|\nabla\vec{n}\|_{L^2(D_2(0))}\\[5mm]
\ds \Vert\nabla\vec{e}_1\Vert_{L^p(D_2(0))}+\Vert\nabla\vec{e}_2\Vert_{L^p(D_2(0))}\;\leq\;C(q)\, \|\nabla\vec{n}\|_{L^p(D_2(0))}\quad.
\end{array}
\ee
We follow step by step the proof of lemma 5.1.4 in \cite{Hel} where $p$ is replacing 2. The main ingredient is the use of lemma~\ref{Lma1} instead of theorem 3.1.7 of \cite{Hel}.
Then the whole argument goes through without any further modification. Once $\vec{e}_i$ satisfying (\ref{conf-fact-lp-1}) is constructed we recall that $\la$ satisfies on $D^2_2(0)$
the equation (\ref{la-equ}). We can then decompose $\la$ into the sum of it's harmonic extension $\la^0$ on $D^2_2$ and the function $\mu$ satisfying (\ref{la-equ}). Applying one more time
lemma~\ref{Lma1} but to $\mu$ this time one obtains
\be
\label{conf-fact-lp-2}
\Vert \nabla\mu\Vert_{L^p(D_2(0))}\;\leq C(q)\ \Vert\nabla\vec{e}_1\Vert_{L^2(D_2(0))}\Vert\nabla\vec{e}_2\Vert_{L^p(D_2(0))}\;\leq\;  C(q)\, \|\nabla\vec{n}\|_{L^2(D_2(0))} \|\nabla\vec{n}\|_{L^p(D_2(0))} \:.
\ee
Standard estimates on harmonic functions give also
\be
\label{conf-fact-lp-3}
\|\nabla\la^0\|_{L^p(D_{3/2}(0))}\le \;C(q)\,\Vert\nabla\la\Vert_{L^{2,\infty}(D_2(0))}
\ee
Combining (\ref{conf-fact-lp-2}) and (\ref{conf-fact-lp-3}) give (\ref{conf-factor-control-Lp}) and lemma~\ref{lm-conf-factor-Lp} is proved.\hfill $\Box$

\subsection{Controlling the $L^p$ norm of $\nabla\vec{n}$.}
We shall now prove the following lemma.
\begin{Lm}
\label{nabla-n-control}
Let $2<q<+\infty$ and consider $\vec{\Phi}\in W^{2,p}_{imm}(D_2(0),{\R}^3)$ for some $2\le p\le q$. If 
\[
\int_{D_2(0)}|\nabla\vec{n}|^2 \ dx^2\le 8\pi/3\quad\mbox{ and }\quad\|\nabla\la\|_{L^{2,\infty}(D_2(0))}=m<+\infty\quad.
\]
where $e^\la:=|\p_{x_1}\vec{\Phi}|=|\p_{x_2}\vec{\Phi}|$. Then, there exists $C(q,m)>0$ such that
\be
\label{nabla-n}
\int_{D_1(0)}|\nabla\vec{n}|^p\ dx^2\le\,C(m,q)\ \lf[ \int_{D_{3/2}(0)}\ e^{2\la}\ dx^2 \rg]^{(p-2)/2}\ \int_{D_{3/2}(0)}|\vec{H}|^p\ e^{2\la}\ dx^2+C(m,q)\quad.
\ee

\hfill $\Box$
\end{Lm}
\noindent{\bf Proof of lemma~\ref{nabla-n-control}.} 
 Denote $\bar{\la}$  the average of $\la$ on $D_{3/2}(0)$.  Using (\ref{conf-factor-control}) together with the fact that
 \[
 \Delta\vec{\Phi}=2\, e^{2\la}\ \vec{H}
 \]
 we obtain
 \be
 \label{conf-fact-1}
 \int_{D_{3/2}(0)}|\Delta\vec{\Phi}|^p\ dx^2\le C(m,q)\ e^{2\bar{\la}\,(p-1)}\ \int_{D_{3/2}(0)}|\vec{H}|^p\ e^{2\la}\ dx^2\quad.
 \ee
 Using classical elliptic estimates we have
 \be
 \label{conf-fact-2}
 \int_{D_1(0)}|\nabla^2\vec{\Phi}|^p\ dx^2\le C(q)\, \int_{D_{3/2}(0)}|\Delta\vec{\Phi}|^p\ dx^2+C(q)\, \lf[ \int_{D_{3/2}(0)}|\nabla\vec{\Phi}|^2\ dx^2 \rg]^{p/2}\quad.
 \ee
The following pointwize estimate holds
 \be
 \label{conf-fact-3}
 |\nabla\vec{n}|(x)=\lf|\nabla\lf(\frac{\p_{x_1}\vec{\Phi}\times\p_{x_2}\vec{\Phi}}{\lf| \p_{x_1}\vec{\Phi}\times\p_{x_2}\vec{\Phi} \rg|}\rg)\rg|(x)\le 4\, |\nabla^2\vec{\Phi}|(x)\ e^{-\la(x)}\quad.
 \ee
This gives obviously using one more time (\ref{conf-factor-control}) 
\be
\label{conf-fact-4}
\int_{D_1(0)}|\nabla\vec{n}|^p\ dx^2\le\,C(m,q)\ e^{-\,p\,\bar{\la}}\ \int_{D_1(0)}|\nabla^2\vec{\Phi}|^p\ dx^2\quad.
\ee
Combining (\ref{conf-fact-2}) and (\ref{conf-fact-4}) we obtain
\be
\label{conf-fact-5}
\begin{array}{l}
\ds\int_{D_1(0)}|\nabla\vec{n}|^p\ dx^2\le\,C(m,q)\ e^{-\,p\,\bar{\la}}\,\int_{D_{3/2}(0)}|\Delta\vec{\Phi}|^p\ dx^2+C(m,q)\ e^{-\,p\,\bar{\la}}\,\lf[ \int_{D_{3/2}(0)}|\nabla\vec{\Phi}|^2\ dx^2 \rg]^{p/2}\\[5mm]\ds\quad\quad\le\, C(m,q)\ e^{-\,p\,\bar{\la}}\,\int_{D_{3/2}(0)}|\Delta\vec{\Phi}|^p\ dx^2+ C(m,q)
\end{array}
\ee
Finally combining (\ref{conf-fact-1}) and (\ref{conf-fact-5}) we obtain
\be
\label{conf-fact-6}
\int_{D_1(0)}|\nabla\vec{n}|^p\ dx^2\le\,C(m,q)\ e^{\bar{\la}\,(p-2)}\ \int_{D_{3/2}(0)}|\vec{H}|^p\ e^{2\la}\ dx^2+C(m,q)\quad.
\ee
This concludes the proof of lemma~\ref{nabla-n-control}.   \hfill $\Box$

\section{Proof of theorem~\ref{th-regularity}.}
\reset

From now on, we take $F(t)$ as in the statement of Theorem \ref{th-regularity}. In this subsection we  establish regularity from a qualitative perspective exclusively.
We shall post-pone quantitative estimates to the next subsection. We treat the case $p>2$ in this subsection. For the most delicate case $p=2$ (already treated in \cite{Riv2})
quantitative estimates as the one derived in the next subsection are needed.

By hypothesis, $\nabla^2\bp$ lies in $L^p$, and thus 
\be\label{com2}
\nabla\bn\,\in\,L^p\:.
\ee
Hence $H$ lie in $L^p$. From hypothesis (a) on $F$, we find that $F(H^2)\text{e}^{2\la}$, and thus $F(H^2)$, are integrable. Per hypothesis (c), it follows that
\be\label{com3}
|HF'|\,\in\,L^{p'}(\di)\:,\qquad\text{where}\:\:\:\quad\dfrac{1}{p'}\,:=\,1-\frac{1}{p}\:,
\ee
so that
\be\label{com1}
H^2F'(H^2)\,\in\, L^p\cdot L^{p'}\,\in\,L^1\:.
\ee
Accordingly, writing (\ref{defL}) in the form
\bes
\nabla^\perp\bL\;=\;-\,\nabla(HF'\bn)+2HF'\nabla\bn+F\nabla\bp
\ees
shows that 
\be\label{regL}
\bL\,\in\,L^{p'}+W^{1,1}\,\subset\,L^{p'}\:.
\ee
We also get from (\ref{defY00}) that $\Delta Y$ is integrable, so that\footnote{$L^{2,\infty}$ is the weak Marcinkiewicz space. Refer to the Appendix for details. That $\nabla Y$ lies in $L^{2,\infty}$ follows from the fact that $\Delta Y$ lies in $L^1$. See Theorem 3.3.6 in \cite{Hel}.}
\be\label{regY}
\nabla Y\,\in\,L^{2,\infty}\:.
\ee
Bringing (\ref{regL}) and (\ref{regY}) into (\ref{defS}) and (\ref{defV}) shows that
\be\label{regSU}
\nabla S\:,\,\nabla\vec{V}\:\in\:L^{p'}+L^{2,\infty}\:\subset\:L^{p'}\:,
\ee
since $p'<2$. \\

%We next recall the following result in integrability by compensation due to Y. Ge in the case $s=2$ \cite{Ge} and to the authors for $s\ne2$ (Lemma IV.2 in \cite{BR1}). 
%\begin{Lma}\label{Lma01}
%Let $D_1(0)$ be the unit disk in ${\R}^2$. 
%Consider the divergence-form problem
%\bes
%\left\{\begin{array}{rclcl}
%\Delta\varphi&=&\nabla^\perp a\cdot \nabla b&,&\text{in}\:\: D_1(0)\\[1ex]
%\varphi&=&0&,&\text{on}\:\:\partial D_1(0)\:,
%\end{array}\right.
%\ees
%where $\nabla a\in L^{2,\infty}(D_1(0))$ and $\nabla b\in L^s(D_1(0))$, for some $s\in(1,\infty)$. There holds
%\bes
%\Vert\nabla\varphi\Vert_{L^s(D_1(0))}\;\leq\;C_{s}\Vert\nabla a\Vert_{L^{2,\infty}(D_1(0))}\Vert\nabla b\Vert_{L^s(D_1(0))}\:,
%\ees
%for some constant $C_s>0$.\\
%Furthermore, when $s\in(2,\infty)$, there holds
%\bes
%\Vert\nabla^2\varphi\Vert_{L^{r}(D_1(0))}+\Vert\nabla\varphi\Vert_{L^s}\;\leq\;C'_{s}\Vert\nabla a\Vert_{L^{2,\infty}}\Vert\nabla b\Vert_{L^s}\:,
%\ees
%for some constant $C'_s>0$ and for all $1\le r<\dfrac{2s}{s+2}$. 
%\end{Lma}
From the Hodge decomposition (\ref{hodgey}), we see that
\be\label{eqv}
\Delta\vec{v}\;=\;\nabla^\perp Y\cdot\nabla\bn\:.
\ee
According to the conditions (\ref{regY}) and (\ref{com2}), since $p>2$, by classical rules on products in Lorentz spaces, we have that $\nabla^\perp Y\cdot\nabla\bn\in L^{2p/(p+2),p}$ and by Lorentz-Sobolev embeddings we obtain
\be\label{regv}
\nabla\vec{v}\,\in\,L^{p}
\qquad\text{and}\qquad\nabla^2\vec{v}\,\in\,L^{r}\quad\text{for all}\:\:\: 1\le r<\dfrac{2p}{p+2} \:,
\ee
which, along with (\ref{regY}) and (\ref{hodgey}), gives us
\be\label{regu}
\nabla\vec{u}\,\in\,L^{2,\infty}\:.
\ee
As $\bR:=\bV-\vec{u}$, we obtain from the latter and (\ref{regSU}) that
\be\label{regSR}
\nabla S\:,\,\nabla\vec{R}\:\in\:L^{p'}\:.
\ee

Next, we differentiate the first order system (\ref{cwsfirst0}) to obtain the second order system in divergence form
\be\label{cws2nd0}
\left\{\begin{array}{rcl}
\Delta\vec{R}&=&\nabla\bn\times\nabla^\perp\bR+\nabla\bn\cdot\nabla^\perp S-\text{div}\big(\bn\times\nabla\vec{v}\big)\\[1ex]
\Delta S&=&-\nabla\bn\cdot\nabla^\perp\bR+\text{div}\big(\bn\cdot\nabla\vec{v}\big) \:.
\end{array}\right.
\ee
We will now call upon another result of integration by compensation  (see \cite{CLMS}). 
\begin{Lma}\label{Lma2}
Let $D:=D_r(x)$ be an arbitrary disk in ${\R}^2$. 
Let $q\in[2,\infty)$, $2\le p\le q$ and let $p'=\dfrac{p}{p-1}$. Consider two functions $a\in W^{1,p}(\di)$ and $b\in W^{1,p'}(\di)$. The solution of the problem
\bes
\left\{\begin{array}{rclcl}
\Delta\varphi&=&\nabla^\perp a\cdot \nabla b&,&\text{in}\:\:\di\\[1ex]
\varphi&=&0&,&\text{on}\:\:\partial\di\:,
\end{array}\right.
\ees
satisfies
\bes
\Vert\nabla^2\varphi\Vert_{L^1(\di)}\;\leq\;C_{q}\Vert\nabla a\Vert_{L^{p}(\di)}\Vert\nabla b\Vert_{L^{p'}(\di)}\:,
\ees
for some constant $C_q>0$ depending only on $q$.\hfill $\Box$
\end{Lma}
With this result at our disposal, we can now use the hypotheses (\ref{com2}), (\ref{regSR}) and simple regularity estimates to derive from the system (\ref{cws2nd0}) that
\be\label{regSR2}
\nabla^2 S\:,\,\nabla^2\vec{R}\:\in\:L^{1}\:.
\ee
In order to weave this information back to the level of curvature, we need a new identity. On one hand, from (\ref{hodgey}), we find
\begin{eqnarray}\label{heu1}
\nabla\vec{V}\times\nabla^\perp\bp&=&\nabla\bR\times\nabla^\perp\bp+\nabla\vec{u}\times\nabla^\perp\bp\;\;=\;\;\nabla\bR\times\nabla^\perp\bp+\big(\nabla\vec{v}-\bn\nabla^\perp Y\big)\times\nabla\bp\nonumber\\[1ex]
&=&\nabla\bR\times\nabla^\perp\bp+\nabla\vec{v}\times\nabla\bp+\nabla Y\cdot\nabla\bp\:.
\end{eqnarray} 
On the other hand, using (\ref{defV0}) and (\ref{defS}), there holds
\begin{eqnarray*}
\nabla\vec{V}\times\nabla^\perp\bp&=&\big(\bL\times\nabla\bp\big)\times\nabla^\perp\bp+HF'(H^2)\nabla\bp\times\nabla^\perp\bp\;\;=\;\;\big(\bL\cdot\nabla^\perp\bp\big)\cdot\nabla\bp-2\text{e}^{2\la}HF'(H^2)\bn\nonumber\\[1ex]
&=&\nabla^\perp S\cdot\nabla\bp+\nabla Y\cdot\nabla\bp-2\text{e}^{2\la}HF'(H^2)\bn\:.
\end{eqnarray*} 
Combining the latter to (\ref{heu1}) yields the identity 
\be\label{back2h}
-\,2\text{e}^{2\la}HF'(H^2)\bn\;=\;\nabla S\cdot\nabla^\perp\bp+\nabla\bR\times\nabla^\perp\bp+\nabla\vec{v}\times\nabla\bp\:.
\ee

Recall that the regularity of the Gauss curvature is tied to that of the Gauss map, namely, $\text{e}^{2\la}K$ is as regular as $\nabla\bn$ is. In particular, the Liouville equation
\be\label{liouville}
-\,\Delta\la\;=\;\text{e}^{2\la}K
\ee
 and the fact that $\bn\in W^{1,p}$ show that $\la\in W^{2,p}\hookrightarrow L^\infty$, and thus that 
\be\label{regla}
\text{e}^{\pm\la}\in W^{2,p}\qquad\text{and}\qquad\nabla\bp\in W^{2,p}\:.
\ee
Combining (\ref{regv}), (\ref{regSR2}), and (\ref{regla}) into (\ref{liouville}) shows that
\bes
HF'(H^2)\,\in\,W^{1,1}\subset L^2\:.
\ees
Since $|HF'(H^2)|\gtrsim|H|^{p-1}$ by hypothesis (c), it follows that
\be\label{regHnew}
 H\,\in\,L^{2(p-1)}\:,
\ee
which, because $p>2$, is an improvement to the original assumption that $H\in L^p$. Recall that the Gauss map satisfies a perturbed harmonic map equation\footnote{cf. Appendix B in \cite{BR2}.}, namely:
\be\label{gauss}
\Delta\bn\;=\;-\,|\nabla\bn|^2\bn-2\,\text{div}\big(H\nabla\bp\big)+4\,\text{e}^{2\la}H^2\bn\:.
\ee
Owing to $\nabla\bn$ being in $L^p$ (for $p>2$) and to (\ref{regHnew}), it follows easily that $\nabla\bn\in L^{2(p-1)}$. We now introduce this information along with (\ref{regSR2}) and (\ref{regv}) into the system (\ref{cws2nd0}) to discover that $\Delta S$ and $\Delta R$ lie in $L^r$ for all $1\le r<\dfrac{2p}{p+2}$. As $p>2$, we obtain the improvement
\be\label{regSR2bis}
\nabla^2 S\:,\,\nabla^2\bR\,\in\,L^r\:,\qquad\forall\:\:\:1\le r<\dfrac{2p}{p+2}\:.
\ee
Just as we did above, we conclude that 
\be\label{oops}
HF'(H^2)\in W^{1,r}\:,\quad\forall\:\:\:r<\dfrac{2p}{p+2}\:,%,\qquad\text{and}\qquad H\,,\,\nabla\bn\in L^s     \:,\quad\forall\:\:\: s<p(p-1)\:.
\ee
and thus
\bes
H\,\in\,L^{q}\qquad\forall\:\:q<{p(p-1)}\:.
\ees
Put into (\ref{gauss}), the latter easily shows that 
\be\label{regn3}
\nabla\bn\,\in\,L^{q}\qquad\forall\:\:q<p(p-1)\:.
\ee
We also have that $F(H^2)$ and $H^2F'(H^2)$ lie in the space $L^s$ for $s<p-1$. Per (\ref{defY00}), it follows that $\nabla Y$ lies in the space $W^{1,s}$, for $s<p-1$. Coupled to (\ref{regn3}) and brought into (\ref{eqv}), this information yields now
\bes
\nabla^2\vec{v}\,\in\,L^a\qquad\forall\:\:\left\{\begin{array}{lcl}a<\dfrac{2p(p-1)}{-p^2+3p+2}&\text{if}&p\le3\\[2ex]
a<p(p-1)&\text{if}&p>3\:.
\end{array}\right.
\ees
This can now be bootstrapped along with (\ref{regn3}) and (\ref{regSR2bis}) back into the system (\ref{cws2nd0}) to obtain that
\bes
\nabla^2S\,,\nabla^2\bR\,\in\,L^s\:,\qquad\forall\:\:\:1\leq s<p-1\:,
\ees
which is an improvement to (\ref{regSR2bis}).\\
The above routine can of course be run again, improving with each step the integrability of the quantities involved. Without much effort, one eventually reaches that 

\be\label{ined}
HF'(H^2)\,\in\,W^{1,s}\:,\qquad\forall\:\:s<\infty\:.
\ee
To transfer this information back to $H$, it suffices to use hypothesis (b), which states that $HF'(H^2)$ is a smooth and invertible function of $H$.  In particular, 
\bes
H\,\in\,W^{1,s}\:,\qquad\forall\:\:s<\infty\:.
\ees
The equation $\Delta\bp=2\text{e}^{2\la}H\bn$ now gives that $\bp\in W^{3,s}$ for all finite $s$. Classical methods easily imply that the immersion $\bp$ is in fact smooth, by repeatedly differentiating all equations and applying standard elliptic estimates. This concludes the proof of theorem~\ref{th-regularity}.

\begin{Rm}\label{casehp}
In the special case when $F(H)=|H|^p$, we have $HF'(H^2)=\frac{p}{2}|H|^{p-1}$, which is not invertible near $H=0$. A glance at the proof reveals however that (\ref{ined}) still holds. In particular, $H$ is a bounded function. 
\end{Rm}

\section{Proof of Theorem \ref{th-epsreg} }\label{enest1}
\reset
The proof is accomplished in several steps. 

\subsection{Improved integrability: from $L^{p',\infty}$ to $L^{p'}$}

We will work within the context of {\it weak immersions with bounded second fundamental form}. That is, we requires $\bp:\Sigma\rightarrow\R^3$ to be Lipschitz and have a non-degenerate pull-back metric $g:=\bp^*g_{\R^3}$. Moreover, we demand that the second fundamental form of the immersion lie in the space $L^2$:
\bes
\int_\Sigma|d\bn|_{g}^2\,d\text{vol}_g\;<\;\infty\:.
\ees
These immersions are well-understood \cite{ParksCity}. In particular, it is known since the work of F. H\'elein, that locally about every point, it is possible to {\it conformally} reparametrize the surface. Since our results will be local in nature, we will thus assume that $\bp$ is conformal, namely, $g=\text{e}^{2\la}\delta$, where $\delta$ denotes the Cartesian metric on $\R^2$, and $\la$ is the conformal parameter. 
Localizing, we assume that there holds\footnote{for simplicity, we have switched to the ``flat" coordinates notation: $|d\bn|_{g}^2\,d\text{vol}_g=|\nabla\bn|^2dx_1dx_2$, where $\nabla:=(\partial_{x_1},\partial_{x_2}).$ }
\be\label{small2}
\int_{D_2(0)}|\nabla\bn|^2\;\leq\;\eps_0
\ee
for some suitably chosen constant $8\pi/3>\eps_0>0$ (to be made precise in the sequel). We will also suppose that
\bes
m\;:=\;\Vert\nabla\la\Vert_{L^{2,\infty}(D_2(0))}\;<\;\infty\:.
\ees
From lemma~\ref{conf-factor} we have
\be
\label{controlfactor}
\lf|\la-\frac{1}{|D_{3/2}(0)|}\int_{D_{3/2}(0)}\la(x)\ dx^2\rg|\le \; C(m)\quad.
\ee
We focus on values of $p$ within the interval $[2,q]$ for some fixed $q>2$. We set $p':=p/(p-1)$. Using lemma~\ref{nabla-n-control} we have
\be\label{1}
\begin{array}{l}
\ds\Vert\nabla\bn\Vert_{L^p(D_1(0))}\;\leq\;C(m,q)\  e^{\ov{\la}\, (1-2/p)}\ \lf[\int_{D_{3/2}(0)}(1+H^2)^{p/2}\ e^{2\la}\ dx^2\rg]^{1/p}+ C(m,q)\\[5mm]
\ds\quad\quad\le C(m,q)\ [e^{\ov{\la}\, (1-2/p)}\ w_p^{1/p}+1]\le C(m,q)\ \sqrt{w_p}\quad.
\end{array}
\ee
where we used $e^{2\,\ov{\la}}\le w_p$. In this subsection, we study critical points $\bp$ of the energy
\bes
\mathcal{W}_p(\Sigma)\;:=\;\int_\Sigma(1+|H|^2)^{p/2}\text{e}^{2\la}dx\qquad\text{with}\qquad w_p\;:=\;\int_{D_{3/2}(0)}(1+|H|^2)^{p/2}\text{e}^{2\la}dx\;<\;\infty\:.
\ees
As we have seen in the first section, critical points of the energy $\mathcal{W}_p$ satisfy 
\be\label{3}
\text{div}\big(-\nabla(HF'\bn\big)+2HF'\nabla\bn+F\nabla\bp \big)\;=\;0\qquad\text{on}\:\:D_{2}(0)\:,
\ee
with $F(H^2)=(1+H^2)^{p/2}$. Note that
\bes
F'\;=\;\dfrac{p}{2}(1+H^2)^{p/2-1}\qquad\text{satisfies}\qquad \Vert F'\Vert_{L^{p/(p-2)}(D_1(0))}\;\leq\;C(q)\ e^{-\,(2-4/p)\,\ov{\la}}\ w_p^{1-2/p}\quad.
\ees
and we have moreover
\be
\label{est-1}
\big\Vert HF'\big\Vert_{L^{p'}(D_1(0))}\le\, e^{-2\, \ov{\la}\,(1-1/p)}\ w_p^{1-1/p}
\ee
Owing to (\ref{1}), we also have
\be\label{3.5}
\big\Vert HF'\nabla\bn\big\Vert_{L^{1}(D_1(0))}\;\leq \ \Vert HF'\big\Vert_{L^{p'}(D_1(0))}\ \Vert\nabla\bn\Vert_{L^p(D_1(0))}\quad .
\ee
We can call upon Lemma \ref{A-1} to integrate (\ref{3}) on the disk $D_1(0)$ and find there exists an element $\bL$ (unique up to the addition of a constant) satisfying
\be\label{4}
\nabla^\perp\bL\;=\;-\nabla\big(HF'\bn\big)+2HF'\nabla\bn+F\nabla\bp\qquad\text{on}\:\: D_1(0)\:,
\ee
and such that
\be
\label{L-control}
\begin{array}{l}
\ds\big\Vert\bL\big\Vert_{L^{p',\infty}(D_1(0))}\;\leq\; C_q\ \big\Vert HF'\big\Vert_{L^{p'}(D_1(0))}+\big\Vert HF'\nabla\bn\big\Vert_{L^{1}(D_1(0))}+\big\Vert F\, e^\la\big\Vert_{L^{1}(D_1(0))}\\[5mm]
\ds\quad\quad\le C(m,q)\ e^{-2\, \ov{\la}\,(1-1/p)}\ w_p^{1-1/p}\ \lf[1+\|\nabla\vec{n}\|_{L^p(D_1(0))}\rg]+\ e^{-\ov{\la}}\ w_p
\end{array}
\ee
where we are choosing $\bL$ to have average 0 on $D_1(0)$. 

Using that $\nabla\bn\cdot\nabla\bp=-2\text{e}^{2\la}H$, we find from (\ref{4})
\be\label{inter1}
\nabla\bp\cdot\nabla^\perp\bL\;=\;2HF'\nabla\bn\cdot\nabla\bp+F|\nabla\bp|^2\;=\;4\text{e}^{2\la}(F-H^2F')\:.
\ee
Let $Y$ be the solution of
\be\label{defY}
\left\{\begin{array}{rclcl}
\Delta Y&=&4\text{e}^{2\la}(F-H^2F')&&\text{in}\:\:D_1(0)\\[1ex]
Y&=&0&&\text{on}\:\:\partial D_1(0)\:.
\end{array}\right.
\ee
Clearly, 
\begin{eqnarray*}
\Vert\Delta Y\Vert_{L^1(D_1(0))}&\leq&C\  \Vert F\, e^{2\la}\Vert_{L^1(\dio)}+C\ \Vert H^2 F'(H)\ e^{2\la}\Vert_{L^{1}(D_1(0))}\le\ C(m)\ w_p\:,
\end{eqnarray*}
so that, owing to a classical linear elliptic result (see e.g. Theorem 3.3.6 in \cite{Hel}):
\be\label{10}
\Vert\nabla Y\Vert_{L^{2,\infty}(D_1(0))}\;\leq\;C(m)\ w_p\:.
\ee

Let next $\bv$ be the unique solution of the following problem:
\be\label{11}
\left\{\begin{array}{rclcl}
\Delta\bv&=&\nabla^\perp Y\cdot\nabla\bn& &\text{in}\:\:D_1(0)\\[1ex]
\bv&=&0&,&\text{on}\:\:\partial D_1(0)\:.
\end{array}\right.
\ee
To obtain estimates on $\bv$, we will use  lemma~\ref{Lma1} in the range $p\in[2,q]$ to the system (\ref{11}) using (\ref{1}) and (\ref{10}). This yields the estimate
\be\label{12}
\begin{array}{l}
\ds\Vert\nabla\bv\Vert_{L^p(D_1(0))}\;\leq\;C_q\ \|\nabla Y\|_{L^{2,\infty}(D_1(0))}\ \|\nabla\vec{n}\|_{L^p(D_1(0))}\ 
\leq\;C(m)\ w_p\ \|\nabla\vec{n}\|_{L^p(D_1(0))}\quad.
\end{array}
\ee
and
\be
\label{12bb}
\ds\Vert\nabla\bv\Vert_{L^2(D_1(0))}\;\leq\;C_q\ \|\nabla Y\|_{L^{2,\infty}(D_1(0))}\ \|\nabla\vec{n}\|_{L^2(D_1(0))}\ 
\leq\;C(q)\ w_p\ \eps_0 \quad.
\ee

\noindent
Note that $\text{div}(\bn\nabla^\perp Y-\nabla\bv)\equiv0$ holds on $D_1(0)$. This equation can be integrated to yield the existence of $\bu$ satisfying
\be\label{13}
\nabla^\perp\vec{u}\;=\;\bn\,\nabla^\perp Y-\nabla\bv\quad.
\ee
Owing to the fact that $L^2\subset L^{2,\infty}$. Combining (\ref{10}) and (\ref{12}) then gives
\be\label{14}
\Vert\nabla\bu\Vert_{L^{2,\infty}(D_1(0))}\;\leq\;C(m, q)\ w_p\quad.
\ee

Let us return to (\ref{inter1}) and (\ref{defY}). Together they give that $\text{div}(\bL\cdot\nabla^\perp\bp-\nabla Y)\equiv0$. Integrating this equation yields the existence of a potential $S$ satisfying
\be\label{15}
\nabla^\perp S\;=\;\bL\cdot\nabla^\perp\bp-\nabla Y\;=\;\bL\cdot\nabla^\perp\bp+\bn\cdot\big(\nabla\bu-\nabla^\perp\bv\big)\:,
\ee
where we have used (\ref{13}).  
Next, from the definition of $\nabla\bL$ given in (\ref{4}) and using the fact that $\bn\times\nabla\bp=-\nabla^\perp\bp$ (and hence that $\nabla\bn\times\nabla\bp=0$), we find
\bes
\text{div}\big(-\bL\times\nabla^\perp\bp  \big)\;=\;\nabla^\perp\bL\times\nabla\bp\;=\;\text{div}\big(HF'\nabla^\perp\bp \big)\:.
\ees
This exact divergence form is also locally integrated to yield a potential function $\bV$ satisfying
\be\label{defV}
\nabla\bV\;=\;\bL\times\nabla\bp+HF'\nabla\bp\:.
\ee
As in the previous section, we set $\bR:=\bV-\vec{u}$ to obtain the system 
\be\label{cws2nd-00}
\left\{\begin{array}{rcl}
\Delta\vec{R}&=&-\nabla^\perp\bn\times\nabla\bR-\nabla^\perp\bn\cdot\nabla S-\text{div}\big(\bn\times\nabla\vec{v}\big)\\[1ex]
\Delta S&=&\nabla^\perp\bn\cdot\nabla\bR+\text{div}\big(\bn\cdot\nabla\vec{v}\big) \:.
\end{array}\right.
\ee
We shall denote
\be
\label{15-a}
S_p:= e^{(1-2/p)\,\ov{\la}}\ S\quad,\quad \vec{R}_p:= e^{(1-2/p)\,\ov{\la}}\ \vec{R}\quad\mbox{ and }\quad \vec{v}_p:=e^{(1-2/p)\,\ov{\la}}\ \vec{v}\quad.
\ee
and by homogeneity we have
\be\label{cws2nd}
\left\{\begin{array}{rcl}
\Delta\vec{R}_p&=&-\nabla^\perp\bn\times\nabla\bR_p-\nabla^\perp\bn\cdot\nabla S_p-\text{div}\big(\bn\times\nabla\vec{v}_p\big)\\[1ex]
\Delta S_p&=&\nabla^\perp\bn\cdot\nabla\bR_p+\text{div}\big(\bn\cdot\nabla\vec{v}_p\big) \:.
\end{array}\right.
\ee
Using the definition of $S$ and $S_p$, combining (\ref{L-control}) and (\ref{10}) we obtain
\be
\label{S_p-control}
\|\nabla S_p\|_{L^{p',\infty}(D_1(0))}\le\  C(m,q)\ w_p^{1-1/p}\ \lf[1+\|\nabla\vec{n}\|_{L^p(D_1(0))}\rg]+C(m)\ e^{(1-2/p)\,\ov{\la}}\ w_p\:
\ee
Recall that $e^{2\, \ov{\la}}\le w_p$, hence, combining (\ref{1}) with (\ref{S_p-control}) we obtain
\be
\label{ov-lambda-control}
\begin{array}{l}
\ds \|\nabla S_p\|_{L^{p',\infty}(D_1(0))}\le\ C(m,q)\ w_p^{1-1/p}\ [e^{(1-2/p)\,\ov{\la}}\ w_p^{1/p}+1]+C(m)\ e^{(1-2/p)\,\ov{\la}}\ w_p\\[5mm]
 \ds\quad\quad\le C(m,q)\ w_p^{1-1/p}\ [\sqrt{w_p}+1]
 \end{array}
 \ee
Combining (\ref{12}) and (\ref{1}) we obtain 
\be
\label{15-b}
\|\nabla\vec{v}_p\|_{L^{p}(D_1(0))}\le C(m)\ e^{(1-2/p)\,\ov{\la}}\ w_p\ \|\nabla\vec{n}\|_{L^p(D_1(0))}\le C(m)\ w_p^{1-1/p}\ [w_p+\sqrt{w_p}]\ 
\ee
and using (\ref{12bb}) we obtain
\be
\label{15-bb-a}
\|\nabla\vec{v}_p\|_{L^{2}(D_1(0))}\le C(m)\ w_p^{1-{1/p}}\ \sqrt{w_p}\quad.
\ee
We also obtain 
\be
\label{15-bb}
\|\nabla\vec{u}_p\|_{L^{2,\infty}(D_1(0))}\le C(m,q)\ w_p^{1-{1/p}}\ \sqrt{w_p}\quad.
\ee
Similarly to $S_p$ we have
\be
\label{15-bc}
\|\nabla\vec{V}_p\|_{L^{p',\infty}(D_1(0))}\le\ C(m,q)\ w_p^{1-1/p}\ [\sqrt{w_p}+1]
\ee
and hence we deduce from (\ref{15-bb}) and (\ref{15-bc})
\be
\label{15-bd}
\|\nabla\vec{R}_p\|_{L^{p',\infty}(D_1(0))}\le\ C(m,q)\ w_p^{1-1/p}\ [\sqrt{w_p}+1]
\ee
Applying lemma~\ref{Lma3} for $q=\infty$ to the system (\ref{cws2nd}) with (\ref{1}), (\ref{15-bb}) and (\ref{15-bd}) we first find
\be\label{add10-l}
\Vert\nabla S_p\Vert_{L^{p',2}(D_{4/5}(0))}+\Vert\nabla \bR_p\Vert_{L^{p',2}(D_{4/5}(0))}\;\leq\; C(m,q)\ w_p^{1-1/p}\ [{w_p}+1]\ \quad.
\ee
and then, applying the same lemma but for $q=2$ we obtain
\be\label{add10}
\Vert\nabla S_p\Vert_{L^{p',1}(D_{3/4}(0))}+\Vert\nabla \bR_p\Vert_{L^{p',1}(D_{3/4}(0))}\;\leq\; C(m,q)\ w_p^{1-1/p}\ [{w_p}+1]\ \quad.
\ee

\subsection{Power decays.}

\paragraph{\underline{The case $p-2$ away from $0$.} }

We assume first that there exists $1>\delta>0$ such that $0<\delta\le p-2\le \delta^{-1}<+\infty$. In this case each step of the bootstrap argument in the proof of theorem~\ref{th-regularity} can be made quantitative with constants depending only on $\delta$.

\paragraph{\underline{The case $p-2$ sufficiently small.} }

This is the most delicate case.
Let $x\in D_{1/2}(0)$ and $r<r_0$ where $r_0$ will be fixed later independent of $x$, $\vec{\Phi}$ and $p\in [2,q]$. We decompose $S_p$ and $\vec{R}_p$ in $D_r(x)$ as follows
\[
S_p=S_p^0+S_p^1\quad\mbox{ where }\lf\{
\begin{array}{l}
\Delta S_p^0=0\quad\mbox{ in }D_r(x)\\[3mm]
S_p^0= S_p\quad\mbox{ on }\p D_r(x)
\end{array}
\rg.\quad\mbox{ and }\quad\vec{R}_p=\vec{R}_p^{\,0}+\vec{R}_p^{\,1}\quad\mbox{ where }\lf\{
\begin{array}{l}
\Delta \vec{R}_p^{\,0}=0\quad\mbox{ in }D_r(x)\\[3mm]
\vec{R}_p^{\,0}= \vec{R}_p\quad\mbox{ on }\p D_r(x)
\end{array}
\rg.
\]
Lemma~\ref{harm-power-decay} implies that for any $t\in (0,1/2)$
\be
\label{powdec-1}
\int_{D_{t\,r}(x)} |\nabla S_p^0|^{p'}(y)+|\nabla \vec{R}_p^{\,0}|^{p'}(y)\ dy^2\le \ C\ t^2\ \int_{D_{r}(x)} |\nabla S_p^0|^{p'}(y)+|\nabla \vec{R}_p^{\,0}|^{p'}(y)\ dy^2
\ee
Using lemma~\ref{Lma1} as well as standard elliptic estimates we have
\be
\label{powdec-2-b}
\begin{array}{l}
\ds\int_{D_{r}(x)} |\nabla S_p^1|^{p'}(y)+|\nabla \vec{R}_p^{\,1}|^{p'}(y)\ dy^2\le C_q\  \eps_0^{p'} \ \int_{D_{r}(x)} |\nabla S_p|^{p'}(y)+|\nabla \vec{R}_p|^{p'}(y)\ dy^2\ \\[3mm]
\ds\quad\quad\quad+\, C_q\ \int_{D_{r}(x)}|\nabla\vec{v}_p|^{p'}(y)\ dy^2
\end{array}
\ee
Combining (\ref{powdec-1}) and (\ref{powdec-2-b}) we obtain
\be
\label{powdec-2}
\begin{array}{l}
\ds\int_{D_{t\,r}(x)} |\nabla S_p|^{p'}(y)+|\nabla \vec{R}_p|^{p'}(y)\ dy^2\le C_q\  (\eps_0^{p'}+t^2) \ \int_{D_{r}(x)} |\nabla S_p|^{p'}(y)+|\nabla \vec{R}_p|^{p'}(y)\ dy^2\ \\[3mm]
\ds\quad\quad\quad+\, C_q\ \int_{D_{r}(x)}|\nabla\vec{v}_p|^{p'}(y)\ dy^2
\end{array}
\ee
We return to (\ref{defY}). Recall that
\bes
F-H^2F'\;=\;F^{(p-2)/p}+\dfrac{1}{p}(p-2)H^2F'\:.
\ees
hence
\be\label{morbihan}
\Delta Y\;=\;4\text{e}^{2\la}F^{(p-2)/p}+\dfrac{4}{p}(p-2)\text{e}^{2\la}H^2F'\qquad\text{on}\:\: D_1(0)\:,\qquad\text{with}\quad Y\big|_{\partial D_1(0)}\,=\,0\:.
\ee
We denote $Y_p:=e^{(1-2/p)\ov{\la}} Y$ and decompose $Y_p=Z_p+W_p$ where
\be
\label{powdec-3}
\lf\{
\begin{array}{l}
\ds\Delta Z_p=4\, e^{(1-2/p)\ov{\la}}\,\text{e}^{2\la}F^{(p-2)/p}\quad\mbox{ in }D_1(0)\\[3mm]
\ds Z_p=0\quad\mbox{ on }\p D_1(0)
\end{array}
\rg.
\ee
We have, for $p/(p-2)>4/3$
\be
\label{powdec-4}
\|\nabla Z_p\|_{L^4(D_1(0))}\le C\ e^{(1+2/p)\ov{\la}}\ \lf[\int_{D_1(0)} F(H^2)\ e^{2\la}\ dx^2\rg]^{(p-2)/p}\le C\ w_p^{1-1/p}\ \sqrt{w_p}\quad.
\ee
This gives in particular that
\be
\label{powdec-5}
\sup_{x\in D_{1/2}(0), \ r<r_0} r^{-1}\ \|\nabla Z_p\|_{L^2(D_r(x))}\le \|\nabla Z_p\|_{L^4(D_1(0))}\le C\ w_p^{1-1/p}\ \sqrt{w_p}\quad.
\ee
Decomposing as above the function $W_p=W^0_p+W_p^1$ on the ball $D_r(x)$ into an harmonic and a trace zero function we obtain, using lemma~\ref{harm-power-decay}, for any $t\in (0,1/2)$
\be
\label{powdec-6}
\|\nabla W_p\|_{L^{2,\infty}(D_{t\,r}(x))}\le\ C\ t\ \|\nabla W_p\|_{L^{2,\infty}(D_{r}(x))}+C(m,q)\, (p-2)\ e^{(1-2/p)\ov{\la}}\ \int_{D_r(x)} H^2\, F'(H^2)\ e^{2\la}\ dy^2
\ee
We write $H^2\, F'(H^2)=H\ H\ F'(H^2)$   and from (\ref{back2h}) we deduce
\be
\label{powdec-7}
\int_{D_r(x)} H^2\, F'(H^2)\ e^{2\la}\ dy^2\le w_p^{1/p}\ \lf[\|\nabla S_p\|_{L^{p'}(D_r(x))}+\|\nabla \vec{R}_p\|_{L^{p'}(D_r(x))}+\|\nabla \vec{v}_p\|_{L^{p'}(D_r(x))}\rg]
\ee
Hence combining (\ref{powdec-5}),  (\ref{powdec-6}) and (\ref{powdec-7}) we obtain
\be
\label{powdec-8}
\begin{array}{l}
\ds\|\nabla Y_p\|_{L^{2,\infty}(D_{t\,r}(x))}\le\ C\ t\ \|\nabla Y_p\|_{L^{2,\infty}(D_{r}(x))}+C\,w_p^{1-1/p}\ \sqrt{w_p}\ r\\[5mm]
\ds\quad\quad+C(m,q)\, (p-2)\  \sqrt{w_p}\ \lf[\|\nabla S_p\|_{L^{p'}(D_r(x))}+\|\nabla \vec{R}_p\|_{L^{p'}(D_r(x))}+\|\nabla \vec{v}_p\|_{L^{p'}(D_r(x))}\rg]
\end{array}
\ee
Decomposing as above the map $\vec{v}_p=\vec{v}^{\,0}_p+\vec{v}_p^{\,1}$ on the ball $D_r(x)$ into an harmonic and a trace zero function we obtain, using lemma~\ref{harm-power-decay}, for any $t\in (0,1/2)$ we obtain, using also lemma~\ref{Lma1},
\be
\label{powdec-9}
\|\nabla \vec{v}_p\|^{p'}_{L^{p'}(D_{t\,r}(x))}\le\ C\ t^2\ \|\nabla \vec{v}_p\|^{p'}_{L^{p'}(D_{r}(x))}+C\ \eps_0^{p'}\ \|\nabla Y_p\|^{p'}_{L^{2,\infty}(D_{r}(x))}
\ee
Combining now (\ref{powdec-2}), (\ref{powdec-8}) as well as (\ref{powdec-9}) for $r$ and for $r/t$ we obtain
\be
\label{powdec-10}
\begin{array}{l}
\ds\int_{D_{t\,r}(x)} |\nabla S_p|^{p'}(y)+|\nabla \vec{R}_p|^{p'}(y)+|\nabla \vec{v}_p|^{p'}(y)\ dy^2+ \|\nabla Y_p\|^{p'}_{L^{2,\infty}(D_{t\,r}(x))}\\[5mm]
\ds\quad\le C(m,q)\  (\eps_0^{p'}+t^2+(p-2)^{p'}\, w_p^{p'/2}) \ \int_{D_{r}(x)} |\nabla S_p|^{p'}(y)+|\nabla \vec{R}_p|^{p'}(y)+|\nabla \vec{v}_p|^{p'}(y)\ dy^2\ \\[5mm]
\ds\quad+\, C_q\ t^{p'}\ \int_{D_{r/t}(x)}|\nabla \vec{v}_p|^{p'}\ dy^2+C\ (\eps_0+t^{p'})\ \|\nabla Y_p\|^{p'}_{L^{2,\infty}(D_{r/t}(x))}+C\,w_p\ w_p^{p'/2}\ r^{p'}
\end{array}
\ee
For $\eps_0$, $t$ and $p-2$ chosen to be small enough we then have
\be
\label{powdec-11}
\begin{array}{l}
\ds\int_{D_{t\,r}(x)} |\nabla S_p|^{p'}(y)+|\nabla \vec{R}_p|^{p'}(y)+|\nabla \vec{v}_p|^{p'}(y)\ dy^2+ \|\nabla Y_p\|^{p'}_{L^{2,\infty}(D_{t\,r}(x))}\\[5mm]
\ds\quad\le 2^{-1}\int_{D_{r/t}(x)} |\nabla S_p|^{p'}(y)+|\nabla \vec{R}_p|^{p'}(y)+|\nabla \vec{v}_p|^{p'}(y)\ dy^2+ \|\nabla Y_p\|^{p'}_{L^{2,\infty}(D_{r/t}(x))}+ C\,w_p\ w_p^{p'/2}\, r^{p'}
\end{array}
\ee
We choose $r_0/t<1/4$ and  by a standard iteration argument, using (\ref{15-bb-a}) and (\ref{add10}) , we deduce the following power decay 
\be
\label{powdec-12}
\begin{array}{l}
\ds\int_{D_{r}(x)} |\nabla S_p|^{p'}(y)+|\nabla \vec{R}_p|^{p'}(y)+|\nabla \vec{v}_p|^{p'}(y)\ dy^2+ \|\nabla Y_p\|^{p'}_{L^{2,\infty}(D_{r}(x))}\\[5mm]
\ds\quad\quad\le\  C(m,q)\, r^\al\ w_p\ [{w_p}+1]^{p'/2}
\end{array}
\ee
where $\al>0$ only depends on $q$ and not on $p\in [2,q]$.

\subsection{Improved integrability and uniform $\eps$-regularity estimate}The obtained power decay, by the mean of (\ref{back2h}), namely
\bes
-2\text{e}^{2\la}HF'\bn\;=\;\nabla S\cdot\nabla^\perp\bp+\nabla\bR\times\nabla^\perp\bp+\nabla\bv\times\nabla\bp\:.
\ees 
gives 
\be\label{H-pointwize}
e^{2\ov{\la}/p'}\ |H|^{p/p'}\le\ C(m,q)\ \lf[|\nabla S_p|+|\nabla\vec{R}_p|+|\nabla\vec{v}_p|\rg]
\ee
Combining (\ref{powdec-12}) and (\ref{H-pointwize}) gives then 
\be
\label{H-dec}
\sup_{x\in D_{1/2}(0)\ ;\ r<r_0}r^{-\al}\int_{D_r(x)} H^p\ e^{2\la}\ dy^2\le C(m,q)\ \ w_p\ [{w_p}+1]^{p'/2}\quad.
\ee
Combining (\ref{H-dec}) with lemma~\ref{nabla-n-control}   we obtain   
\be
\label{HFprimenablan}
\sup_{x\in D_{1/2}(0)\ ;\ r<r_0}r^{-\al/p'}\int_{D_r(x)} |H\, F'(H)\ \nabla\vec{n}|\ dy^2\le \ C(m,q)\ e^{-2\la/p'}\ w_p^{1/p'}\ (\sqrt{w_p}+1)\ \sqrt{w_p}\quad.
\ee
moreover we have
\be
\label{FnablaPhi}
\sup_{x\in D_{1/2}(0)\ ;\ r<r_0}r^{-\al}\int_{D_r(x)} F(H)\ |\nabla\vec{\Phi}|\ dy^2\le C(m,q)\ e^{-\ov{\la}}\ w_p\ [{w_p}+1]^{p'/2}\quad.
\ee
Inserting these power decays in the r-h-s of the Willmore PDE
\be\label{powdec-13}
\Delta(HF'\bn)=\text{div}\big(2HF'\nabla\bn+F\nabla\bp \big)\qquad\text{on}\:\:D_{3/4}(0)\:,
\ee
gives, by the mean of Adams estimates \cite{Ad2}, a higher integrability to $HF'\bn$ of the form
\[
\|H\, F'\|_{L^{r}(D_{1/4}(0))}\le C(m,q)\ \lf[e^{-2\la/p'}\ w_p^{1/p'}\ (\sqrt{w_p}+1)\ \sqrt{w_p}+e^{-\ov{\la}}\ w_p\ [{w_p}+1]^{p'/2}\rg]
\]
for some $r>p'$. We can now bootstrap this estimate in (\ref{powdec-13}) to reach the desired $\epsilon$-regularity statement 
\[
e^{2\la(0)} H^p(0)\le C(m,q)\ w_p\ \lf[(1+\sqrt{w_p})^{p'}+{w_p}^{p'/2}\ [1+w_p]^{(p')^2/2}\rg]\quad.
\]

\section{Proof of Theorem \ref{epsreg2} }\label{enest2}
\reset
\subsection{Rewriting the Euler Lagrange Equation of ${\mathcal E}_p$ and preliminary estimates.}
We place ourselves in the same context as in Section \ref{enest1}. We assume $2\le p\le 3$. We consider a local conformal parametrization $\bp:D_2(0)\rightarrow\R^3$ with conformal parameter $\la$. We suppose that
\be\label{small3}
\int_{D_2(0)}|\nabla\bn|^2\;\leq\;\eps_0^2\:,
\ee
for some suitably chosen constant $\eps_0>0$ (to be made precise in the sequel). We will also suppose that
\bes
m\;:=\;\Vert\nabla\la\Vert_{L^{2,\infty}(D_2(0))}\;<\;\infty\:.
\ees
Following lemma~\ref{conf-factor} we have
\be
\label{conf-factor-control-bis}
\|\nabla\la\|_{L^2(D_{3/2}(0))}+\lf|\la-\frac{1}{|D_{3/2}(0)|}\int_{D_{3/2}(0)}\la(x)\ dx^2\rg|\le \;C\,\Vert\nabla\la\Vert_{L^{2,\infty}(D_2(0))}+\;  C\, \|\nabla\vec{n}\|^2_{L^2(D_2(0))} \quad.
\ee
H\'elein's construction of the moving frame (see \cite{Hel} lemma 5.1.4) gives not only a Coulomb frame $(\vec{e}_1,\vec{e}_2)$ satisfying (\ref{coulomb}) but 
it also satisfies 
\be
\label{coulomb-p}
\|\nabla\vec{e}_1\|_{L^p(D_2(0))}+\|\nabla\vec{e}_2\|_{L^p(D_2(0))}\le\ C\ \|\nabla\vec{n}\|_{L^p(D_2(0))}
\ee
where $C$ is independent of $p\in [2,3]$. Indeed, it is proved that there exists $\nu\in W^{1,2}_0(D_2)$ such that
\[
\nabla\nu=(\vec{e}_1,\nabla^\perp\vec{e}_2)
\]
Hence $\nu$ is the solution of
\[
\lf\{
\begin{array}{l}
\ds \Delta\nu=(\nabla\vec{e}_1;\nabla^\perp\vec{e}_2)\quad\mbox{ in }D_2(0)\\[3mm]
\ds \nu=0\quad\mbox{ on }\partial D_2(0)
\end{array}
\rg.
\]
Using Lemma~\ref{Lma1} we have the a-priori estimate (using (\ref{coulomb}))
\be
\label{apriori-coulomb-p}
\|\nabla\nu\|_{L^p(D_2(0))}\le\ C\ \|\nabla \vec{e}_1\|_{L^2(D_2(0))}\ \|\nabla \vec{e}_2\|_{L^p(D_2(0))}\le \, C\ \varepsilon_0\ \|\nabla \vec{e}_2\|_{L^p(D_2(0))}
\ee
where $C$ is independent of $p\in [2,3]$. We have also the pointwize identity
\be
\label{apriori-coulomb-p-2-a}
|\nabla\vec{e}_1|^2+|\nabla\vec{e}_2|^2= 2\,|\nabla\nu|^2+|\nabla\vec{n}|^2\quad.
\ee
Hence we deduce the a-priori inequality
\be
\label{apriori-coulomb-p-2}
\|\nabla\vec{e}_1\|_{L^p(D_2(0))}+\|\nabla\vec{e}_2\|_{L^p(D_2(0))}\le\ C\ \varepsilon_0\ \|\nabla \vec{e}_2\|_{L^p(D_2(0))}+\|\nabla\vec{n}\|_{L^p(D_2(0))}
\ee
and (\ref{coulomb-p}) follows by a density argument as in \cite{Hel}. With (\ref{coulomb-p}) at hand, using one more time the full strength of lemma~\ref{Lma1}, we can adapt step by step the proof of lemma~\ref{conf-factor} to deduce
\be
\label{conf-factor-control-p}
\|\nabla\la\|_{L^p(D_{3/2}(0))}\le \;C\,\Vert\nabla\la\Vert_{L^{2,\infty}(D_2(0))}+\;  C\, \|\nabla\vec{n}\|_{L^p(D_2(0))}\quad,
\ee
where $C$ is \underbar{independent} of $p\in [2,3]$.

\medskip

As in the proof of Lemma~\ref{nabla-n-control} we have
\be
\label{hessien-p}
\int_{D_{6/5}(0)}|\nabla^2\vec{\Phi}|^p\ dx^2\le C(m)\ e^{2\ov{\la}(p-1)}\ \int_{D_2(0)}|\vec{H}|^p\ e^{2\la}\ dx^2+C\ \lf[\int_{D_{3/2}(0)}|\nabla\vec{\Phi}|^2\ dx^2\rg]^{p/2}\quad.
\ee
We set $p':=p/(p-1)$, and in this section, we study critical points $\bp$ of the energy
\bes
\mathcal{E}_p(\Sigma)\;:=\;\int_{\Sigma}(1+|A|^2)^{p/2}\,d\text{vol}_g\qquad\text{with}\qquad e_p\;:=\;\int_{D_2(0)}(1+|A|^{2})^{p/2}\text{e}^{2\la}\,dx\;<\;\infty\:.
\ees
For notational convenience, we will set
\bes
F\;=\;(1+|A|^2)^{p/2}\qquad\text{and}\qquad F'\;=\;\dfrac{p}{2}\,(1+|A|^2)^{p/2-1}\:.
\ees
where $|A|^2=e^{-2\la}\,|\nabla\vec{n}|^2$. With these notations we have 
\be
\label{hessien-p-p}
\|\nabla^2\vec{\Phi}\|_{L^p(D_{1}(0))}\le \ C(m)\ e^{2\,\ov{\la}/p'}\ e_p^{1/p}+ C\ e^{\ov{\la}}
\ee
For future reference, note that
\be\label{add1000}
\Vert F'(|A|^2)\,|\nabla\vec{n}|^2\Vert_{L^1(D_1(0))}\;\leq\; e_p
\ee
and
\be\label{add1001}
\Vert F'(|A|^2)\, |\nabla\vec{n}|\Vert_{L^{p'}(D_1(0))}\le\; e^{\ov{\la}\,(1-2/p')}\ e_p^{1/p'}
\ee
As we have seen in Section \ref{vary}, critical points of $\mathcal{E}_p$ satisfy (\ref{hoopla}):
\bes
d^{*_g}\Big[Fd\bp-2F'(d\bn\stackrel{.}{\otimes}d\bn)\res_gd\bp+2\big(d^{*_g}(F'd\bn)\cdot d\bp\big)\bn \Big]\;=\;0\:.
\ees
In conformal coordinates, this expression reads
\be\label{divee}
\text{div}\Big[F\nabla\bp-2\,\text{e}^{-2\la}F'\sum_{j=1}^{2}(\nabla\bn\cdot\partial_{x_j}\bn)\partial_{x_j}\bp+2\text{e}^{-2\la}\bn\big(\text{div}\big(F'\nabla\bn\big)\cdot\nabla\bp\big)  \Big]\;=\;0\:.
\ee
We have
\bes
\Big|F\nabla\bp-2\,\text{e}^{-2\la}F'\sum_{j=1}^{2}(\nabla\bn\cdot\partial_{x_j}\bn)\partial_{x_j}\bp\Big|\;\leq\;\text{e}^{\la}\,F+2\, e^{-\la}\,|\nabla\bn|^2\,F'\:,
\ees 
so that owing to (\ref{add1000}), 
\be\label{add30}
\Big\Vert F\nabla\bp-2\,\text{e}^{-2\la}F'\sum_{j=1}^{2}(\nabla\bn\cdot\partial_{x_j}\bn)\partial_{x_j}\bp\Big\Vert_{L^1(D_1(0))}\le\,C(m)\, e^{-\ov{\la}}\ e_p\quad.
\ee where we have used that $p'\leq2$. 
We also have
\be\label{add31}
\big\Vert \text{div}(F'\nabla\bn)\big\Vert_{W^{-1,p'}(D_{1}(0))}\;\leq\;\Vert F'\nabla\bn\Vert_{L^{p'}(D_1(0))}\;\leq\;e^{\ov{\la}\,(1-2/p')}\ e_p^{1/p'}\:,
\ee
where we have used (\ref{add1001}). From (\ref{conf-factor-control-p}) we have
\be
\label{algeb-1}
\|e^\la\|_{L^\infty\cap W^{1,p}(D_1(0))}\le\ C(m)\ \lf[e^{\ov{\la}}+e^{2\,\ov{\la}/p'} e_p^{1/p}\rg] \quad,
\ee
from (\ref{hessien-p-p})
\be
\label{algeb-2}
\|\nabla\vec{\Phi}\|_{L^\infty\cap W^{1,p}(D_1(0))}\le\ C(m)\ \lf[e^{\ov{\la}}+e^{2\,\ov{\la}/p'} e_p^{1/p}\rg] \quad,
\ee
and
\be
\label{algeb-3}
\|\vec{n}\|_{L^\infty\cap W^{1,p}(D_1(0))}\le \ 1+e^{-2\ov{\la}/p}\ e_p^{1/p}
\ee
where $C(m)$ is \underbar{independent} of $p$. We have also
\be
\label{algeb-4}
\|e^{-2\la}\ \vec{n}\otimes \nabla\vec{\Phi}\|_{L^\infty\cap W^{1,p}(D_1(0))}\le \, \ C(m)\ \lf[e^{-\ov{\la}}+e^{-\,2\,\ov{\la}/p} e_p^{1/p}\rg] \quad.
\ee
Recall that there exists a constant $C$ independent of $p\in[2,3]$ such that for all $f\in C_0^\infty$
\[
\|f\|_{L^\infty(D_1(0))}\le C\ \|\nabla f\|_{L^{p,1}(D_1(0))}\quad.
\]
Hence, combining (\ref{add1001}) and (\ref{algeb-4})
\be
\label{algeb-5}
\sup_{\|f\|_{W^{1,(p,1)}_0(D_1(0))}\le1}\ \int_{D_1(0)}f(x)\ \text{e}^{-2\la}\bn\big(\text{div}(F'\nabla\bn)\cdot\nabla\bp\big)\ dx^2\le C(m)\ \lf[e^{-2\,\ov{\la}/p'}\,e_p^{1/p'}+e^{-\,\ov{\la}}\ e_p\rg]\ 
\ee
Hence
\be\label{add32}
\big\Vert \text{e}^{-2\la}\bn\big(\text{div}(F'\nabla\bn)\cdot\nabla\bp\big)\big\Vert_{W^{-1,(p',\infty)}(D_{1}(0))}\;\leq\; C(m)\ \lf[e^{-2\,\ov{\la}/p'}\,e_p^{1/p'}+e^{-\,\ov{\la}}\ e_p\rg]\ 
\ee
Let $\vec{L}$ be the distribution (unique up to a constant) such that
\be
\label{algeb-6}
\nabla^\perp\bL\;=\;F\nabla\bp-2\text{e}^{-2\la}F'\sum_{j=1}^{2}(\nabla\bn\cdot\partial_{x_j}\bn)\partial_{x_j}\bp+2\text{e}^{-2\la}\bn\big(\text{div}(F'\nabla\bn)\cdot\nabla\bp\big) 
\ee
Combining (\ref{add30}) and (\ref{add32}) we obtain
\be
\label{algeb-7}
\|\nabla^\perp\bL\|_{L^1(D_1(0))+W^{-1,(p',\infty)}(D_1(0))}\le\ C(m)\ \lf[e^{-2\,\ov{\la}/p'}\,e_p^{1/p'}+e^{-\,\ov{\la}}\ e_p\rg]\quad,
\ee
from which we deduce
 \be\label{estimL}
\Vert\bL\Vert_{L^{p',\infty}(D_{1}(0))}\;\leq\;C(m)\ \lf[e^{-2\,\ov{\la}/p'}\,e_p^{1/p'}+e^{-\,\ov{\la}}\ e_p\rg]
\ee
where we normalized $\vec{L}$ to have average 0 on $D_1(0)$. We compute
\begin{eqnarray}\label{aux1}
\nabla\bp\cdot\nabla^\perp\vec{L}&=&2\text{e}^{2\la}F-2F'|\nabla\bn|^2\nonumber\\[1ex]
&=&\text{e}^{2\la}\big(1+\text{e}^{-2\la}|\nabla\bn|^2\big)^{p/2-1}\big(2+(2-p)\text{e}^{-2\la}|\nabla\bn|^2\big)\:.
\end{eqnarray}
where we have used that $|\nabla\bn|^2=\text{e}^{2\la}A^2$. 
Let us next introduce
\be\label{defY2}
\left\{\begin{array}{rclcl}
-\Delta Y&=&2\text{e}^{2\la}F-2F'|\nabla\bn|^2&,&\text{in}\:\:D_1(0)\\[1ex]
Y&=&0&,&\text{on}\:\:\partial D_1(0)\:.
\end{array}\right.
\ee
According to (\ref{aux1}), we thus have an exact divergence form which may be locally integrated to yield a new potential $S$ satisfying
\be\label{defSX}
\nabla S\;=\;\,\bL\cdot\nabla\bp+\nabla^\perp Y\:.
\ee
We next compute
\begin{eqnarray}\label{aux2}
\nabla\bp\times\nabla^\perp\vec{L}&=&-2\,\text{e}^{-2\la}\big(\nabla\bp\cdot\text{div}(F' \nabla\bn)\big)\cdot(\bn\times\nabla\bp)\nonumber\\[1ex]
&=&\:-2\,\text{e}^{-2\la}\sum_{j=1}^{2}\big(\partial_{x_j}\bp\cdot\text{div}(F' \nabla\bn)\big)(\bn\times\partial_{x_j}\bp)\nonumber\\[1ex]
&=&\:-2\,\bn\times\text{div}\big(F' \nabla\bn\big)\;\;=\;\;-\,\text{div}\Big[2F'\,\bn\times\nabla\bn \Big]\nonumber\\[1ex]
&=&\:-\,\text{div}\Big[2F'\nabla^\perp\bn+4\,HF'\nabla^\perp\bp   \Big]\:,
\end{eqnarray}
where we have used that $\nabla\bp\times\nabla\bn=0\,$ and that $\nabla\bp\times\bn=\nabla^\perp\bp$. We have also used the elementary identity
\be\label{aux3}
\nabla^\perp\bn\;=\;\bn\times\nabla\bn-2H\nabla^\perp\bp\:.
\ee
Equivalently, (\ref{aux2}) reads
\be\label{aux10}
\text{div}\Big[\vec{L}\times\nabla^\perp\bp+2F'\,\nabla^\perp\bn+4\,HF'\nabla^\perp\bp \Big]\;=\;0\:.
\ee
Integrating  this identity yields a potential $\bV$ satisfying
\be\label{defR1}
\nabla^\perp\bV  \;=\;\vec{L}\times\nabla^\perp\bp+2F'\nabla^\perp\bn+4HF'\nabla^\perp\bp \:.
\ee
Observe that
\be
\label{vs-1}
\vec{n}\cdot\nabla\vec{V}=\vec{n}\cdot (\vec{L}\times\nabla\vec{\Phi})=\vec{L}\cdot(\nabla\vec{\Phi}\times \vec{n})=\vec{L}\cdot\nabla^\perp\vec{\Phi}=\nabla^\perp S+\nabla Y
\ee
We decompose
\be
\label{vs-2}
2\, F'(|A|^2)\, \nabla\vec{n}=-\nabla\vec{\sigma}+\nabla^\perp\vec{\tau}
\ee
where $\vec{\tau}$ is zero on $\p D_1(0)$ and
\be
\label{vs-3}
\vec{n}\ \nabla Y=\nabla \vec{u}-\nabla^{\perp}\vec{w}
\ee
where $\vec{w}$ is equal to zero on $\p D_1(0)$. We have
\be
\label{vs-4}
\begin{array}{l}
\ds\nabla(\vec{V}+\vec{\sigma})=\vec{L}\times\nabla\vec{\Phi}+4\ H\, F'(|A|^2)\, \nabla\vec{\Phi}+\nabla^\perp\vec{\tau}\\[3mm]
\ds\quad=-(\vec{L}\cdot\vec{n})\ \nabla^\perp\vec{\Phi}+\vec{L}\cdot\nabla^\perp\vec{\Phi}\,\vec{n}+4\ H\, F'(|A|^2)\, \nabla\vec{\Phi}+\nabla^\perp\vec{\tau}\\[3mm]
\ds\quad=-(\vec{L}\cdot\vec{n})\ \nabla^\perp\vec{\Phi}-\vec{n}\,\nabla^\perp S-\vec{n}\,\nabla Y+4\ H\, F'(|A|^2)\, \nabla\vec{\Phi}+\nabla^\perp\vec{\tau}\\[3mm]
\ds\quad=-(\vec{L}\cdot\vec{n})\ \nabla^\perp\vec{\Phi}+4\ H\, F'(|A|^2)\, \nabla\vec{\Phi}-\vec{n}\,\nabla^\perp S-\nabla\vec{u}+\nabla^\perp(\vec{\tau}+\vec{w})
\end{array}
\ee
Denote
\[
\vec{R}:=\vec{V}+\vec{\sigma}+\vec{u}\quad\mbox{ and }\quad\vec{v}:=\vec{\tau}+\vec{w}
\]
We have the following
\begin{Prop}
\label{pr-III.1}
With the previous notations the following equation hold
\be
\label{sys-n}
\vec{n}\times\nabla\vec{R}=-\nabla^\perp\vec{R}+\vec{n}\,\nabla S-\nabla\vec{v}+\vec{n}\times\nabla^\perp\vec{v}\quad.
\ee
\hfill$\Box$
\end{Prop}
\noindent{\bf Proof of proposition~\ref{pr-III.1}.} Observe that
\be
\label{vs-5}
\lf\{
\begin{array}{l}
\ds\vec{n}\times\nabla^\perp\vec{\Phi}=\nabla\vec{\Phi}\\[3mm]
\ds\vec{n}\times\nabla\vec{\Phi}=-\nabla^\perp\vec{\Phi}
\end{array}
\rg.
\ee
Using these identities we deduce in one hand
\be
\label{vs-5aa}
\vec{n}\times \nabla\vec{R}=-(\vec{L}\cdot\vec{n})\ \nabla\vec{\Phi}-4\ H\, F'(|A|^2)\, \nabla^\perp\vec{\Phi}+\vec{n}\times\nabla^\perp\vec{v}\quad,
\ee
and in the other hand
\be
\label{vs-5bb}
\nabla^\perp\vec{R}=(\vec{L}\cdot\vec{n})\ \nabla\vec{\Phi}+4\ H\, F'(|A|^2)\, \nabla^\perp\vec{\Phi}+\vec{n}\,\nabla S-\nabla\vec{v}\quad.
\ee
Combining the two previous lines gives the proposition.\hfill $\Box$

\subsection{improved integrability : from $L^{p',\infty}$ to $L^{p'}$.}
Let $r_0<1/4$ to be fixed later. We proceed to the following renormalization
\[
Y_p:=e^{(1-2/p)\,\ov{\la}}\, Y\quad,\quad S_p:=e^{(1-2/p)\,\ov{\la}}\, S\quad,\quad\bR_p:=e^{(1-2/p)\,\ov{\la}}\, \bR\quad\mbox{and}\quad\vec{v}_p:=e^{(1-2/p)\,\ov{\la}}\, \vec{v}\quad.
\]
Using (\ref{defY2}) we have 
\be
\label{vs-6}
\|\nabla Y_p\|_{L^{2,\infty}(D_1(0))}\le C\ e_p^{1-1/p}\ \sqrt{e_p}\quad.
\ee
Using (\ref{estimL}), (\ref{defSX}) and (\ref{vs-6}) we have
\be
\label{vs-7}
\|\nabla S_p\|_{L^{p',\infty}(D_1(0))}\le C(m)\ e_p^{1/p'}+C\ e_p^{1-1/p}\ \sqrt{e_p}\quad.
\ee
Using (\ref{add1001}), (\ref{estimL}), (\ref{defR1}) and the fact that $|H|^2\le e^{-2\,\la}\,|\nabla\vec{n}|^2$ we have
\be
\label{vs-8}
\|\nabla \vec{V}_p\|_{L^{p',\infty}(D_1(0))}\le C(m)\ e_p^{1/p'}+C\ e_p^{1-1/p}\ \sqrt{e_p}\quad.
\ee
From (\ref{vs-6}) we deduce using standard elliptic estimates
\be
\label{vs-9}
\|\nabla \vec{u}_p\|_{L^{2,\infty}(D_1(0))}+\|\nabla \vec{w}_p\|_{L^{2,\infty}(D_1(0))}\le\ C\ \|\nabla Y_p\|_{L^{2,\infty}(D_1(0))}\le C\ e_p^{1-1/p}\ \sqrt{e_p}\quad.
\ee
where we used the fact that $\vec{w}_p$ satisfies the following system
\be
\label{vs-10}
\lf\{
\begin{array}{l}
\ds -\Delta \vec{w}_p=\nabla^\perp Y_p\nabla\vec{n}\quad\mbox{ in }D_1(0)\\[5mm]
\ds \vec{w}_p=0\quad\mbox{ on }D_1(0)
\end{array}
\rg.
\ee
Using lemma~\ref{Lma1} we deduce 
\be
\label{vs-11}
\|\nabla \vec{w}_p\|_{L^{2}(D_1(0))}\le\ C\ \|\nabla Y_p\|_{L^{2,\infty}(D_1(0))}\ \|\nabla\vec{n}\|_{L^{2}(D_1(0))}\le C\ \varepsilon_0\ e_p^{1-1/p}\ \sqrt{e_p}\
\ee
and from (\ref{add1001}) and (\ref{vs-2}) we have using again standard elliptic estimates
\be
\label{vs-12}
\|\nabla \vec{\sigma}_p\|_{L^{p'}(D_1(0))}+\|\nabla \vec{\tau}_p\|_{L^{p'}(D_1(0))}\le\ C\ e^{\ov{\la}\,(1-2/p)}\ \Vert F'(|A|^2)\, |\nabla\vec{n}|\Vert_{L^{p'}(D_1(0))}\le\; C(m)\ e_p^{1/p'}
\ee
where the constant $C$ is independent on $p\in [2,3]$. Hence we deduce from (\ref{vs-8}), (\ref{vs-9}) and (\ref{vs-12})
\be
\label{vs-12-b}
\|\nabla \vec{R}_p\|_{L^{p',\infty}(D_1(0))}\le C(m)\ e_p^{1/p'}+C\ e_p^{1-1/p}\ \sqrt{e_p}\quad.
\ee
Combining (\ref{vs-11}) and (\ref{vs-12}) we deduce
\be
\label{vs-13}
\|\nabla \vec{v}_p\|_{L^{p'}(D_1(0))}\le\; C(m)\ e_p^{1/p'}+C\ \varepsilon_0\ e_p^{1-1/p}\ \sqrt{e_p}\quad.
\ee

From proposition~\ref{pr-III.1} we have that $(S_p,\vec{R}_p)$ satisfies the following system
\be
\label{vs-14}
\lf\{
\begin{array}{l}
\ds\Delta\vec{R}_p=\nabla S_p\cdot\nabla^\perp\vec{n}+\nabla\vec{R}_p\times\nabla^\perp\vec{n}+\mbox{div}(\vec{n}\times\nabla\vec{v}_p)\\[5mm]
\ds\Delta S_p=\nabla \vec{R}_p\cdot\nabla^\perp\vec{n}+\mbox{div}(\vec{n}\cdot\nabla\vec{v}_p)
\end{array}
\rg.
\ee 
Applying lemma~\ref{Lma3} for $q=\infty$ to the system (\ref{vs-14}) with (\ref{vs-7}), (\ref{vs-12-b}) and (\ref{vs-13})
\be\label{vs-15}
\Vert\nabla S_p\Vert_{L^{p',2}(D_{4/5}(0))}+\Vert\nabla \bR_p\Vert_{L^{p',2}(D_{4/5}(0))}\;\leq\; C(m)\ e_p^{1-1/p}\ [{e_p}+1]\ \quad.
\ee
and then, applying the same lemma but for $q=2$ we obtain
\be\label{vs-16}
\Vert\nabla S_p\Vert_{L^{p',1}(D_{3/4}(0))}+\Vert\nabla \bR_p\Vert_{L^{p',1}(D_{3/4}(0))}\;\leq\; C(m)\ e_p^{1-1/p}\ [{e_p}+1]\ \quad.
\ee

\subsection{Power decays.}
Consider now an arbitrary point $x\in D_{1/2}(0)$ and a radius $0<r<r_0$ where $r_0<1/4$ will be fixed later. Exactly as in the proof of (\ref{powdec-2}), by exploiting the system
(\ref{vs-14}) one first established that for $t\in(0,1)$
\be
\label{vs-16-b}
\begin{array}{l}
\ds\int_{D_{t\,r}(x)} |\nabla S_p|^{p'}(y)+|\nabla \vec{R}_p|^{p'}(y)\ dy^2\le C\  (\eps_0^{p'}+t^2) \ \int_{D_{r}(x)} |\nabla S_p|^{p'}(y)+|\nabla \vec{R}_p|^{p'}(y)\ dy^2\ \\[3mm]
\ds\quad\quad\quad+\, C\ \int_{D_{r}(x)}|\nabla\vec{v}_p|^{p'}(y)\ dy^2
\end{array}
\ee
where $C$ is independent of $p\in [2,3]$.

\medskip

As in the previous section for the Lagrangian ${\mathcal W}_p$, we  proceed to the following decomposition  $Y_p:=Z_p+W_p$ where
\[
\lf\{
\begin{array}{l}
\ds-\Delta Z_p=2\, e^{(1-2/p)\,\ov{\la}}\,\text{e}^{2\la}\big(1+\text{e}^{-2\la}|\nabla\bn|^2\big)^{p/2-1}\quad\mbox{ in }D_1(0)\\[5mm]
\ds Z_p=0\quad\mbox{ on }\p D_1(0)
\end{array}
\rg.
\]
We have for $p/(p-2)>4/3$
\be
\label{zp}
\begin{array}{l}
\ds\|\nabla Z_p\|_{L^4(D_1(0))}\le C(m)\, e^{(1+2/p)\,\ov{\la}}\, \lf[\int_{D_1(0)}(1+e^{-2\la}|\nabla\vec{n}|^2)^{p/2}\ e^{2\la}\ dx^2\rg]^{(p-2)/p}\\[5mm]
\ds\quad\quad\quad\quad\le C(m)\, e_p^{1-1/p}\,\sqrt{e_p}\quad.
\end{array}
\ee
This implies in particular
\be
\label{zp-decay}
\sup_{x\in D_{1/2}(0)\, ;\, r<r_0}r^{-1}\ \|\nabla Z_p\|_{L^2(D_r(x))}\le C(m)\, e_p^{1-1/p}\,\sqrt{e_p}\quad.
\ee
By decomposing $W_p=W_p^0+W_p^1$ as in the previous section into the sum of an harmonic and a trace free part on $\p D_r(x)$ we obtain a similar estimate as (\ref{powdec-6})
\be
\label{vs-17}
\|\nabla W_p\|_{L^{2,\infty}(D_{t\,r}(x))}\le\ C\ t\ \|\nabla W_p\|_{L^{2,\infty}(D_{r}(x))}+C(m)\, (p-2)\ e^{(1-2/p)\ov{\la}}\ \int_{D_r(x)}\big(1+\text{e}^{-2\la}|\nabla\bn|^2\big)^{p/2}\ e^{2\la}\ dy^2\quad.
\ee
Combining (\ref{zp-decay}) and (\ref{vs-17}) gives then
\be
\label{vs-18}
\begin{array}{l}
\ds\|\nabla Y_p\|_{L^{2,\infty}(D_{t\,r}(x))}\le\ C\ t\ \|\nabla Y_p\|_{L^{2,\infty}(D_{r}(x))}\\[5mm]
\ds\quad\quad+C(m)\, (p-2)\ e^{(1-2/p)\ov{\la}}\ \int_{D_r(x)}\big(1+\text{e}^{-2\la}|\nabla\bn|^2\big)^{p/2}\ e^{2\la}\ dy^2+C\ r\quad.
\end{array}
\ee

\medskip

Taking the dot of (\ref{vs-5bb}) with $\nabla^\perp\vec{\Phi}$ we obtain
\be
\label{vs-19}
8\ e^{2\la}\ H\ F'(|A|^2)=\nabla\vec{\Phi}\cdot\nabla\vec{R}-\nabla\vec{\Phi}\cdot\nabla^\perp\vec{v}\quad,
\ee
from which we deduce
\be
\label{vs-20}
\int_{D_r(x)}|H|^p\ e^{2\la}\ dy^2\le C(m)\ \int_{D_r(x)} |\nabla \vec{R}_p|^{p'}+|\nabla\vec{v}_p|^{p'}\ dy^2\quad.
\ee

\medskip

Recall the structural equation
\be
\label{vs-22}
-\nabla\vec{n}=\vec{n}\times\nabla^\perp\vec{n}+\,2\, H\,\nabla\vec{\Phi}\quad.
\ee
Taking the divergence gives
\be
\label{vs-22-a}
-\Delta\vec{n}=\nabla\vec{n}\times\nabla^\perp\vec{n}+2\ \mbox{div}(H\,\nabla\vec{\Phi})
\ee
Decomposing again $\vec{n}$ into the sum of an harmonic part and a trace free part in $D_r(x)$ and making use of the lemma~\ref{Lma1} as well as lemma~\ref{harm-power-decay}
we obtain for any $t\in (0,1)$
\be
\label{vs-23}
\begin{array}{l}
\ds\int_{D_{t\,r}(x)}|\nabla\vec{n}|^p\ e^{-\la\, p}\ e^{2\la}\le C(m)\ (t^2+\varepsilon^2)\ \int_{D_r(x)}|\nabla\vec{n}|^p\ e^{-\la\, p}\ e^{2\la}+ C(m)\ \int_{D_r(x)}|H|^p\ e^{2\la}\ dy^2
\end{array}
\ee
Combining (\ref{vs-20}) and (\ref{vs-23}) gives then
\be
\label{vs-24}
\begin{array}{l}
\ds\int_{D_{t\,r}(x)}|\nabla\vec{n}|^p\ e^{-\la\, p}\ e^{2\la}\ dy^2\le C(m)\ (t^2+\varepsilon^2)\  \int_{D_r(x)}|\nabla\vec{n}|^p\ e^{-\la\, p}\ e^{2\la}\\[5mm]
\ds\quad\quad+ C(m)\ \int_{D_r(x)} |\nabla \vec{R}_p|^{p'}+|\nabla\vec{v}_p|^{p'}\ dy^2\quad.
\end{array}
\ee
Observe that since $\la-\ov{\la}$ is bounded by a constant $C(m)$ one has
\be
\label{vs-25}
\ds\int_{D_{t\,r}(x)} e^{2\la}\ dy^2\le C(m)\ t^2\ \int_{D_r(x)} e^{2\la}\ dy^2
\ee
Hence, combining (\ref{vs-24}) and (\ref{vs-25}) we have
\be
\label{vs-26}
\begin{array}{l}
\ds\int_{D_{t\,r}(x)}(1+|\nabla\vec{n}|^2\ e^{-2\,\la})^{p/2}\ e^{2\la}\ dy^2\le C(m)\ (t^2+\varepsilon^2)\ \ \int_{D_r(x)}(1+|\nabla\vec{n}|^2\ e^{-2\,\la})^{p/2}\ e^{2\la}\\[5mm]
\ds\quad\quad+ C(m)\ \int_{D_r(x)} |\nabla \vec{R}_p|^{p'}+|\nabla\vec{v}_p|^{p'}\ dy^2\quad.
\end{array}
\ee

\medskip

Now it remains to establish some decay of the $L^{p'}$ norm of $\nabla\vec{v}_p$. Recall $\nabla\vec{v}_p:=\nabla\vec{w}_p+\nabla\vec{\tau}_p$. For the $\nabla\vec{w}_p$
the argument follows the lines similar to what we did in the previous section for $\vec{v}_p$ : Using equation (\ref{vs-10}) together with  lemma~\ref{harm-power-decay}, for any $t\in (0,1)$ we obtain, using also lemma~\ref{Lma1}, exactly as we did for proving (\ref{powdec-9})
\be
\label{vs-27}
\|\nabla \vec{w}_p\|^{p'}_{L^{p'}(D_{t\,r}(x))}\le\ C\ t^2\ \|\nabla \vec{w}_p\|^{p'}_{L^{p'}(D_{r}(x))}+C\ \eps^{p'}_0\ \|\nabla Y_p\|^{p'}_{L^{2,\infty}(D_{r}(x))}
\ee
where $C$ is independent of $p\in [2,3]$.

In order to establish some interesting decay for the $L^{p'}$ norm of $\vec{\tau}_p$ we have to use some fundamental lemma in $p-$harmonic theory which is the only
new ingredient in this section in comparison with the previous one. We decompose $\vec{\tau}_p:=\vec{\tau}_p^{\,0}+\vec{\tau}_p^{\,1}$ where $\vec{\tau}^{\,0}_p$ is harmonic
in $D_r(x)$ and equal to $\vec{\tau}_p$ on $\p D_r(x)$. Hence, using (\ref{vs-2}), $\vec{\tau}_p^{\,1}$ satisfies
\be
\label{vs-28}
\lf\{
\begin{array}{l}
\ds-\Delta\vec{\tau}_p^{\,1}=p\,e^{\frac{2-p}{p'}\,\ov{\la}}\ \mbox{div}\lf(  \lf(e^{2\ov{\la}}+e^{-2\, (\la-\ov{\la})}|\nabla^\perp\vec{n}|^2\rg)^{p/2-1}\ \nabla^\perp\vec{n}  \rg)   \quad\quad\mbox{ in } D_r(x) \\[5mm]
\ds \vec{\tau}_p^{\,1}=0\quad\quad\mbox{ on }\p D_r(x)
\end{array}
\rg.
\ee
Using (\ref{A-16}) we have for any $p-2<\al<1/2$
\be
\label{vs-29}
\|\nabla\vec{\tau}_p^{\,1}\|_{L^{p/(p-1)}(D^2_r(x))}\le  C\ e^{ \frac{2-p}{p'}\,\ov{\la}}\  \frac{(p-2)}{\al}\ \lf(\|\nabla\vec{n}\|^{p-1}_{L^p(D^2_r(x))}+(r^{2/p}\,e^{\ov{\la}})^{(p-2)+\al}\, \|\nabla\vec{n}\|^{1-\al}_{L^p(D^2_r(x))}\rg)\quad,
\ee
Observe that $1-(1-\al)/(p-1)=(p-2+\al)/(p-1)$ thus $a^{p-2+\al}\ b^{(1-\al)/(p-1)}\le C_\al\ [a^{p-1}+b]$ where $C_\al$ is uniformly bounded for $\al<1/2$ bounded from below away
from 0. We choose $\al=1/4$ and we have
\be
\label{vs-30}
\|\nabla\vec{\tau}_p^{\,1}\|_{L^{p/(p-1)}(D^2_r(x))}\le  C(m)\ e^{ \frac{2-p}{p'}\,\ov{\la}}\  {(p-2)}\ \lf(\|\nabla\vec{n}\|^{p-1}_{L^p(D^2_r(x))}+(r^{2}\,e^{p\,\ov{\la}})^{(p-1)/p}\rg)\quad.
\ee
Which gives for any $t<1$
\be
\label{vs-31}
\begin{array}{l}
\ds\int_{D_{t\,r}(x))}|\nabla\vec{\tau}_p|^{p'}(y)\ dy^2\le C\ t^2\ \int_{D_{r}(x))}|\nabla\vec{\tau}_p|^{p'}(y)\ dy^2\\[5mm]
\ds\quad\quad+\, C\ \lf[{p-2} \rg]^{p'}\ e^{(2-p)\,\ov{\la}}\ \int_{D_{r}(x)}|\nabla\vec{n}|^p\ dy^2+\, C\ \lf(p-2\rg)^{p'}\ r^2\, e^{2\,\ov{\la}}\\[5mm]
\ds\quad\quad\le C\ t^2\ \int_{D_{r}(x))}|\nabla\vec{\tau}_p|^{p'}(y)\ dy^2+C\ \lf[{p-2} \rg]^{p'}\ \int_{D_r(x)}(1+|\nabla\vec{n}|^2\ e^{-2\,\la})^{p/2}\ e^{2\la}\ dy^2
\end{array}
\ee
Combining (\ref{vs-16-b}), (\ref{vs-18}), (\ref{vs-26}), (\ref{vs-27}) and (\ref{vs-31}) we obtain
\be
\label{vs-32}
\begin{array}{l}
\ds\int_{D_{t\,r}(x)} |\nabla S_p|^{p'}+|\nabla \vec{R}_p|^{p'}+|\nabla\vec{w}_p|^{p'}+|\nabla\vec{\tau}_p|^{p'}+(1+|\nabla\vec{n}|^2\ e^{-2\,\la})^{p/2}\ e^{2\la}\ dy^2+\|\nabla Y_p\|^{p'}_{L^{2,\infty}(D_{t\,r})}\\[5mm]
\ds\le \delta \lf[\int_{D_{\frac{r}{t}}(x)} |\nabla S_p|^{p'}+|\nabla \vec{R}_p|^{p'}+|\nabla\vec{w}_p|^{p'}+|\nabla\vec{\tau}_p|^{p'}+(1+|\nabla\vec{n}|^2\ e^{-2\,\la})^{p/2}\ e^{2\la}\ dy^2+\|\nabla Y_p\|^{p'}_{L^{2,\infty}(D_{\frac{r}{t}})}\rg]\\[3mm]
\ds\quad +C\ r^{p'}
\end{array}
\ee
where $\delta= C(m)\  (\eps_0^{p'}+t^2+t^{p'}+(p-2)^{p'})$. Taking $\eps_0$, $t$ and $p-2$ small enough we can choose $\delta=1/2$ and we deduce, as in the previous section, the following power decay
\be
\label{vs-32-a}
\begin{array}{l}
\ds\sup_{x\in D_{1/2}(0)\,;\, r<1/4}r^{-\al}\ \int_{D^2_r(x)}|\nabla S_p|^{p'}+|\nabla \vec{R}_p|^{p'}+|\nabla\vec{w}_p|^{p'}+|\nabla\vec{\tau}_p|^{p'}+(1+|\nabla\vec{n}|^2\ e^{-2\,\la})^{p/2}\ e^{2\la}\ dy^2\\[3mm]
\ds\quad\quad\quad\quad\le C(m)\ e_p\ [1+e_p^{p'}]
\end{array}
\ee
for some $\al>0$ independent of $p\simeq 2$, where we also used (\ref{vs-13}) and (\ref{vs-16}).

\subsection{Improved integrability and uniform $\eps$-regularity estimate}

First of all we establish improved integrability for the following quantities : $\vec{n}\cdot\nabla\vec{R}_p$ , $\vec{n}\cdot\nabla\vec{v}_p$, $\vec{n}\times\nabla\vec{R}_p+\nabla\vec{v}_p$, $\nabla S_p$.
We deduce from (\ref{algeb-6})
\be
\label{vs-33}
\vec{n}\cdot\nabla^\perp\vec{R}_p=\nabla S_p-\vec{n}\cdot\nabla\vec{v}_p\quad
\ee
This gives in particular
\be
\label{vs-34}
\lf\{
\begin{array}{l}
\ds \mbox{div}\lf(\vec{n}\cdot\nabla\vec{R}_p\rg)=\nabla\vec{n}\cdot\nabla^\perp\vec{v}_p\\[3mm]
\ds \mbox{curl}\lf(\vec{n}\cdot\nabla\vec{R}_p\rg)=-\, \nabla\vec{n}\cdot\nabla^\perp\vec{R}_p
\end{array}
\rg.
\ee
Observe that
\be
\label{vs-35}
\int_{D_2^2}(1+|\nabla\vec{n}|^2\ e^{-2\,\la})^{p/2}\ e^{2\la}\ dy^2\le e_p\quad\Rightarrow\quad e^{2\ov{\la}}\le C(m)\ e_p\quad\mbox{and}\quad\int_{D_1^2}|\nabla\vec{n}|^p\ dx^2\le\, C(m)\ e_p^{p/2}\quad.
\ee
Combining (\ref{vs-32-a}) and (\ref{vs-34}) and (\ref{vs-35}) gives
\be
\label{vs-36}
\sup_{x\in D_{1/2}(0)\,;\, r<1/4}r^{-\al/p'}\int_{D^2_r(x)}\lf| \mbox{div}\lf(\vec{n}\cdot\nabla\vec{R}_p\rg)\rg|+\lf|\mbox{curl}\lf(\vec{n}\cdot\nabla\vec{R}_p\rg)\rg|\le C(m)\  e_p^{1/p'}\ [1+e_p]\ e_p^{1/2}
\ee
Adams elliptic estimates from \cite{Ad2} gives then that for any $q<\frac{2-\al/p'}{1-\al/p'}$
\be
\label{vs-37}
\|\vec{n}\cdot\nabla\vec{R}_p\|_{L^q(D^2_{1/4})}\le C(m)\  e_p^{1/p'}\ [1+e_p]\ [1+e_p^{1/2}]\quad.
\ee
From (\ref{vs-2}) and (\ref{vs-3}) we have
\be
\label{vs-38}
\vec{n}\cdot\nabla^\perp\vec{v}_p=-\nabla Y_p+\vec{n}\cdot\nabla\lf(\vec{\sigma}_p+\vec{u}_p\rg)
\ee
This gives in particular
\be
\label{vs-39}
\lf\{
\begin{array}{l}
\ds \mbox{div}\lf(\vec{n}\cdot\nabla\vec{v}_p\rg)=-\nabla\vec{n}\cdot\nabla^\perp({\sigma}_p+\vec{u}_p)\\[3mm]
\ds \mbox{curl}\lf(\vec{n}\cdot\nabla\vec{v}_p\rg)=-\, \nabla\vec{n}\cdot\nabla^\perp\vec{v}_p
\end{array}
\rg.
\ee
and arguing exactly as for $\vec{n}\cdot\nabla\vec{R}_p$ we obtain that for any $q<\frac{2-\al/p'}{1-\al/p'}$
\be
\label{vs-40}
\|\vec{n}\cdot\nabla\vec{v}_p\|_{L^q(D^2_{1/4})}\le C(m)\  e_p^{1/p'}\ [1+e_p]\ [1+e_p^{1/2}]\quad.
\ee
Combining (\ref{vs-33}), (\ref{vs-37}) and (\ref{vs-40}) we obtain that for any $q<\frac{2-\al/p'}{1-\al/p'}$
\be
\label{vs-40-b}
\|\nabla S_p\|_{L^q(D^2_{1/4})}\le C(m)\  e_p^{1/p'}\ [1+e_p]\ [1+e_p^{1/2}]\quad.
\ee
Using now (\ref{sys-n}), we have
\be
\label{vs-41}
\lf\{
\begin{array}{l}
\ds \mbox{div}\lf(\vec{n}\times\nabla\vec{R}_p+\nabla\vec{v}_p\rg)=-\nabla\vec{n}\times\nabla^\perp\vec{v}_p+\mbox{ div}(\nabla S_p\, \vec{n})\\[3mm]
\ds \mbox{curl}\lf(\vec{n}\times\nabla\vec{R}_p+\nabla\vec{v}_p\rg)=-\, \nabla\vec{n}\cdot\nabla^\perp\vec{R}_p
\end{array}
\rg.
\ee
Combinining (\ref{vs-41}) and (\ref{vs-40-b}) we obtain as before that for any $q<\frac{2-\al/p'}{1-\al/p'}$
\be
\label{vs-42}
\|\vec{n}\times\nabla\vec{R}_p+\nabla\vec{v}_p\|_{L^q(D^2_{1/4})}\le C(m)\  e_p^{1/p'}\ [1+e_p]\ [1+e_p^{1/2}]\quad.
\ee
and combining (\ref{vs-37}) and (\ref{vs-42}) we obtain
 that for any $q<\frac{2-\al/p'}{1-\al/p'}$
\be
\label{vs-43}
\|\vec{n}\times\nabla\vec{v}_p-\nabla\vec{R}_p\|_{L^q(D^2_{1/4})}\le C(m)\  e_p^{1/p'}\ [1+e_p]\ [1+e_p^{1/2}]\quad.
\ee
Combining (\ref{vs-36}) and (\ref{vs-10}) we obtain
\be
\label{vs-44}
\sup_{x\in D_{1/2}(0)\,;\, r<1/4}r^{-\al/p'}\int_{D^2_r(x)}\lf|\Delta\vec{w}_p\rg|\ dy^2\le C(m)\  e_p^{1/p'}\ [1+e_p]\ e_p^{1/2}
\ee
and we deduce as before that for any $q<\frac{2-\al/p'}{1-\al/p'}$
\be
\label{vs-45}
\|\nabla\vec{w}_p\|_{L^q(D^2_{1/4})}\le C(m)\  e_p^{1/p'}\ [1+e_p]\ [1+e_p^{1/2}]\quad.
\ee
From (\ref{vs-22-a}), using again the decomposition of $\vec{n}=\vec{n}^{\,0}+\vec{n}^{\,1}$ into an harmonic part and a trace free part, we deduce by standard elliptic estimates that,
for any $x\in D_{1/8}(0)$ and any $r<1/16$,
\be
\label{vs-46}
\begin{array}{l}
\ds\int_{D_{r/2}(x)}|\nabla\vec{n}|^p\le C\ \lf[\int_{D_r(x)}|\nabla\vec{n}|^{2\,p/3}\ dy^2\rg]^{3/2}+C\ \int_{D_r(x)}H^p\ e^{p\la}\ dy^2
\end{array}
\ee
From which we deduce (still for a constant $C(m)>0$ independent of $r$, $x$ and $3\ge p\ge 2$ but possibly depending on $m$)
\be
\label{vs-47}
\ds\dashint_{D_{r/2}(x)}|\nabla\vec{n}|^p\ e^{(2-p)\,\la}\ dy^2\le C(m)\ \lf[\dashint_{D_r(x)}|\nabla\vec{n}|^{2\,p/3}\ e^{(4-2p)/3}dy^2\rg]^{3/2}+C(m)\ \dashint_{D_r(x)}H^p\ e^{2\la}\ dy^2
\ee
Using (\ref{vs-20}), we have also
\be
\label{vs-48}
\dashint_{D_r(x)}|H|^p\ e^{2\la}\ dy^2\le C(m)\ \dashint_{D_r(x)} \lf[|\nabla \vec{R}_p-\vec{n}\times\nabla\vec{v}_p|^{p'}+|\nabla\vec{w}_p|^{p'}+|\nabla\vec{\tau}_p|^{p'}\rg]\ dy^2\quad.
\ee
Decomposing again $\vec{\tau}:=\vec{\tau}^{\,0}+\vec{\tau}^{\,1}$ into an harmonic part and a trace free part, using (\ref{vs-30}), we deduce from (\ref{vs-30}) by standard elliptic estimates that,
for any $x\in D_{1/8}(0)$ and any $r<1/16$,
\be
\label{vs-49}
\begin{array}{l}
\ds\dashint_{D_r(x)} |\nabla\vec{\tau}_p|^{p'}\ dy^2\le C(m)\ \lf[\dashint_{D_{2r}(x)}|\nabla\vec{\tau}_p|^{3/2}\ dy^2\rg]^{2\,p'/3}\\[5mm]
\ds\quad\quad+C(m)\ \lf[{p-2} \rg]^{p'}\ \dashint_{D_{2r}(x)}(1+|\nabla\vec{n}|^2\ e^{-2\,\la})^{p/2}\ e^{2\la}\ dy^2
\end{array}
\ee
Combining (\ref{vs-47}), (\ref{vs-48}) and (\ref{vs-49}) we obtain
\be
\label{vs-50}
\begin{array}{l}
\ds\dashint_{D_{r/2}(x)}\lf[(1+|\nabla\vec{n}|^2\ e^{-2\,\la})^{p/2}\ e^{2\la}+ |\nabla\vec{\tau}_p|^{p'}\rg]\ \ dy^2\\[5mm]
\ds\quad \quad\le\,C(m)\ \lf[\dashint_{D_{2r}(x)}  \lf[\lf[(1+|\nabla\vec{n}|^2\ e^{-2\,\la}\rg)^{p}\ e^{2\la}+ |\nabla\vec{\tau}_p|^{p'}\rg]^{2/3}\ dy^2\rg]^{3/2}\\[5mm]
\ds\quad\quad\ +\,C(m)\ \dashint_{D_{2\,r}(x)} \lf[|\nabla \vec{R}_p-\vec{n}\times\nabla\vec{v}_p|^{p'}+|\nabla\vec{w}_p|^{p'}\rg]\ dy^2\\[5mm]
\ds\quad\quad+\,C(m)\ \lf[{p-2} \rg]^{p'}\ \dashint_{D_{2r}(x)}(1+|\nabla\vec{n}|^2\ e^{-2\,\la})^{p/2}\ e^{2\la}\ dy^2
\end{array}
\ee
Using Gehring type lemma given by proposition 1.1 of \cite{Gia} as well as estimates (\ref{vs-43}) and (\ref{vs-45}), we obtain for $p-2$ small enough the existence of $q>\max\{p,2+\beta\}$ where $\beta$ is independent of $p$ larger or equal than 2 and sufficiently close to 2 such that
\be
\label{vs-51}
\lf\||\nabla\vec{n}|\,e^{\la \,(2-p)/p}\rg\|_{L^q(D_{1/8}(0))}\le C(m)\ \  e_p^{1/p'}\ [1+e_p]\ [1+e_p^{1/2}]\quad.
\ee
This implies theorem~\ref{epsreg2}. \hfill $\Box$

\renewcommand{\theequation}{A.\arabic{equation}}
\renewcommand{\theTh}{A.\arabic{Th}}
\renewcommand{\theProp}{A.\arabic{Prop}}
\renewcommand{\theLma}{A.\arabic{Lma}}
\renewcommand{\theCo}{A.\arabic{Co}}
\renewcommand{\theRm}{A.\arabic{Rm}}
\renewcommand{\theequation}{A.\arabic{equation}}
\setcounter{equation}{0} 
\reset
\appendix
\section{Appendix}

The following lemma is well known but we give a proof for the reader's convenience.
\begin{Lma}
\label{harm-power-decay}
Let $f$ be an harmonic function on the unit disc $D_1(0)$ , then for any $p\in [1,2]$ and any $t\in (0,1)$ the following inequality holds
\be
\label{dec-harm-1}
\int_{D_t(0)}| f|^p\ dx^2\le C\ t^2\ \int_{D_1(0)}| f|^p\ dx^2\quad,
\ee 
where $C>0$ is universal. Similarly,  the following inequality holds for any $t\in(0,1)$
\be
\label{dec-harm-2}
\| f\|_{L^{2,\infty}(D_t(0))}\le C\ t\ \| f\|_{L^{2,\infty}(D_1(0))}\quad.
\ee
\hfill $\Box$
\end{Lma}
\noindent{\bf Proof of lemma~\ref{harm-power-decay}}

Since $f$ is harmonic we have
\[
\Delta |f|^2=2\,|\nabla f|^2\ge0
\]
This implies that
\[
\int_{\p D_r(0)}\p_r|f|^2\ge 0
\]
from which we deduce that
\[
\frac{d}{d r}\lf[\frac{1}{r^2}\int_{D_r(0)}|f|^2\ dx^2\rg]=\frac{d}{d r}\lf[\int_{D_1(0)}|f(r\,x)|^2\ dx^2\rg]=r\,\int_{0}^1\, ds\int_{\p D_s(0)}\p_r|f|^2(rx)\ge 0
\]
Let $t<1/2$ we have using H\"older for $p\in[1,2]$
\[
\int_{D_t(0)}|f|^p\ dx^2\le \pi^{1-p/2}\,t^{2-p}\ \lf[  \int_{D_t(0)}|f|^2(x)\ dx^2 \rg]^{p/2}\le \pi^{1-p/2}\,t^{2-p}\ \lf[ 4\,t^2\, \int_{D_{1/2}(0)}|f|^2(x)\ dx^2 \rg]^{p/2}
\]
Since $f$ is harmonic, classical elliptic estimates give
\[
\lf[\int_{D_{1/2}(0)}|f|^2(x)\ dx^2 \rg]^{p/2}\le C\ \int_{D_{1}(0)}|f|^p(x)\ dx^2 \quad.
\]
where $C$ is a universal constant. Combining the two previous identities gives (\ref{dec-harm-1}). The inequality (\ref{dec-harm-2}) is derived in a similar way.\hfill $\Box$

\medskip

We now prove the following integrability by compensation result due to Y. Ge in the case $s=2$ \cite{Ge} and to the authors for $s\ne2$ (Lemma IV.2 in \cite{BR1}). 
\begin{Lma}\label{Lma1}
Let $D$ be a disk.
Consider the divergence-form problem
\bes
\left\{\begin{array}{rclcl}
\Delta\varphi&=&\nabla^\perp a\cdot \nabla b&,&\text{in}\:\:\di\\[1ex]
\varphi&=&0&,&\text{on}\:\:\partial\di\:,
\end{array}\right.
\ees
where $\nabla a\in L^{2,\infty}(\di)$ and $\nabla b\in L^p(\di)$, for some $p\in(1,\infty)$. There holds
\be\label{u1}
\Vert\nabla\varphi\Vert_{L^p(\di)}\;\leq\;C_p\Vert\nabla a\Vert_{L^{2,\infty}(\di)}\Vert\nabla b\Vert_{L^p(\di)}\:,
\ee
for some constant $C_p>0$ depending only on $p$ which satisfies the estimate
\be\label{u3}
C_p\;\leq\;C(p_0,p_1)\qquad\forall\:\:1<p_0\leq p\leq p_1<\infty\:,
\ee
for some constant $C(p_0,p_1)$ independent of $p$. \\
%Furthermore, when $s\in(2,\infty)$, there holds
%\be\label{u2}
%\Vert\nabla^2\varphi\Vert_{L^{r}(\di)}+\Vert\nabla\varphi\Vert_{L^s(\di)}\;\leq\;C'_s\Vert\nabla a\Vert_{L^{2,\infty}(\di)}\Vert\nabla b\Vert_{L^s(\di)}\:,
%\ee
%for some constant $C'_s>0$ and for all $1\le r<\dfrac{2s}{s+2}$.
\end{Lma}
{\bf Proof.} The inequality (\ref{u1}) is proved in details in the aforementioned references. 
%As for (\ref{u2}), it is a consequence of $p>2$ and the well-known H\"older inequality for Lorentz spaces. Namely, $L^{2,\infty}\cdot L^p\subset L^{r}$ continuously, for all $1\leq r<2p/(p+2)$. It then suffices to apply the standard Calderon-Zygmund estimate to $\Delta\varphi$ to deduce (\ref{u2}). \\
The estimate (\ref{u3}) can be reached using the Marcinkiewicz interpolation theorem. Let $1<p_0\leq p\leq p_1<\infty$. As the operator $\nabla b\mapsto\Vert\nabla a\Vert^{-1}_{L^{2,\infty}}\nabla\varphi$ is strongly continuous from $L^{(p_0+1)/2}$ to itself and from $L^{p_1+1}$ to itself with norms $M_{(p_0+1)/2}$ and $M_{p_1+1}$ respectively, then it is also strongly continuous from $L^p$ to $L^p$ with norm 
\bes
p^{1/p}\bigg(\dfrac{M_{(p_0+1)/2}^{(p_0+1)/2}}{p-(p_0+1)/2}+\dfrac{M_{p_1+1}^{p_1+1}}{p_1+1-p}\bigg)\;\leq\;\dfrac{3}{2}\bigg(\dfrac{2}{p_0-1}\,M_{(p_0+1)/2}^{(p_0+1)/2}+M_{p_1+1}^{p_1+1}\bigg)\:,
\ees
which depends only on $p_0$ and $p_1$, as desired. \\[-1.5ex]

$\hfill\Box$ \\

We will need a Wente-type result. 

\begin{Lma}\label{Lemma10}
Let $q\in(2,\infty)$ and $p\in[2,q]$. 
Consider the divergence-form problem
\bes
\left\{\begin{array}{rclcl}
\Delta\varphi&=&\nabla^\perp a\cdot \nabla b&,&\text{in}\:\:\di\\[1ex]
\varphi&=&0&,&\text{on}\:\:\partial\di\:,
\end{array}\right.
\ees
where $\nabla a\in L^{p}(\di)$ and $\nabla b\in L^{p',\infty}(\di)$, for some $s\in(1,\infty)$ with $p':=p/(p-1)$. There holds
\be\label{trag}
%\Vert\nabla\varphi\Vert_{L^{2,\infty}(\di)}\;\leq\;
\Vert\nabla\varphi\Vert_{L^{2,p}(\di)}\;\leq\;C_{q}\Vert\nabla a\Vert_{L^{p}(\di)}\Vert\nabla b\Vert_{L^{p',\infty}(\di)}\:,
\ee
for some constant $C_q>0$ depending only on $q$.
\end{Lma}
{\bf Proof.} The estimate (\ref{trag}) holds for $p=2$, as shown by Y. Ge in \cite{Ge}. Let now $p>2$ so that $p'<2$. By a continuous embedding of L. Tartar (see Lemma 31.2 in \cite{Tar1}) and the H\"older inequality for Lorentz spaces, we find:
\bes
b\nabla^\perp a\;\in\;W^{1,(p',\infty)}\cdot L^{p}\;\subset\; L^{\frac{2p}{p-2},\infty}\cdot L^{p,p}\;\subset\;L^{2,p}\:.
\ees
Since $\Delta\varphi=\text{div}(b\nabla^\perp a)$, the estimate (\ref{trag}) follows accordingly for $p>2$. \\
Now, if $p<2$ so that $p'>2$, the Sobolev-Lorentz embedding theorem along with the H\"older inequality for Lorentz spaces, we find
\bes
a\nabla^\perp b\;\in\;W^{1,s}\cdot L^{p',\infty}\;\subset\; L^{\frac{2p}{2-p},p}\cdot L^{p',\infty}\;\subset\;L^{2,p}\:,
\ees
which again proves (\ref{trag}) for $p<2$, since $\Delta\varphi=-\text{div}(a\nabla^\perp b)$. \\

Consider the operator
\bes
B\,:\,(\nabla a,\nabla b)\:\longmapsto\:\nabla\varphi\:.
\ees
As we have seen above, $B$ maps $L^{q'}\times L^{q,\infty}$ continuously into $L^{2,q'}$ with norm $M_0(q)$ ; and it maps continuously $L^{q}\times L^{q',\infty}$ into $L^{2,q}$ with norm $M_1(q)$. 
The bilinear interpolation result of J.-L. Lions and J. Peetre given in Lemma 28.3 of \cite{Tar1} implies now that $B$ maps continuously the spaces\footnote{See \cite{Tar1} for the notation and basic results about interpolation spaces.}
\bes
\big(L^{q'},L^{q}\big)_{\theta,p}\times \big(L^{q,\infty},L^{q',\infty}\big)_{\theta,\infty}\;\equiv\;L^{(\frac{1-\theta}{q'}+\frac{\theta}{q})^{-1},p}  \times L^{(\frac{1-\theta}{q}+\frac{\theta}{q'})^{-1},\infty}
\ees
into 
\bes
\big(L^{2,q'},L^{2,q}\big)_{\theta,p}\;\equiv\;L^{2,p}\:,
\ees
for all $0\leq\theta\leq1$ and $p^{-1}=(1-\theta)/q'+\theta/q$. Equivalently, $B$ maps continuously
\bes
L^{p}\times L^{p',\infty}\quad\text{into}\quad L^{2,p}\qquad\text{for}\quad p\in[q',q]
\ees
with a norm $C_q=C_q(M_0(q),M_1(q),q)$ depending only on $q$, owing to the nature of interpolation. \\[-1.5ex]

$\hfill\Box$ \\

\begin{Lma}\label{Lma3}
Let $D$ be a disk.
Consider the divergence-form problem
\bes
\left\{\begin{array}{rclcl}
\Delta\varphi&=&\nabla^\perp a\cdot \nabla b&,&\text{in}\:\:\di\\[1ex]
\varphi&=&0&,&\text{on}\:\:\partial\di\:,
\end{array}\right.
\ees
where $\nabla a\in L^{2}(\di)$ and $\nabla b\in L^{p,q}(\di)$, for some $p\in(1,\infty)$ and $2\le q\le\infty$. There holds
\be\label{uuu1}
\Vert\nabla\varphi\Vert_{L^{p,q^\ast}(\di)}\;\leq\;C_{p}\Vert\nabla a\Vert_{L^{2}(\di)}\Vert\nabla b\Vert_{L^{p,q}(\di)}\:,
\ee
where $1/q^\ast=1/2+1/q$ for some constant $C_{p}>0$ depending only on $p$ which satisfies the estimate
\be\label{uuu3}
C_{p}\;\leq\;C(p_0,p_1)\qquad\forall\:\:1<p_0\leq p\leq p_1<\infty\:,
\ee
for some constant $C(p_0,p_1)$ independent of $p$.
\hfill $\Box$

%Furthermore, when $s\in(2,\infty)$, there holds
%\be\label{u2}
%\Vert\nabla^2\varphi\Vert_{L^{r}(\di)}+\Vert\nabla\varphi\Vert_{L^s(\di)}\;\leq\;C'_s\Vert\nabla a\Vert_{L^{2,\infty}(\di)}\Vert\nabla b\Vert_{L^s(\di)}\:,
%\ee
%for some constant $C'_s>0$ and for all $1\le r<\dfrac{2s}{s+2}$.
\end{Lma}
{\bf Proof of lemma~\ref{Lma3}.} It suffices to consider the case when $\infty>p_1>2$ and $1<p_0<2$. According to the H\"older inequality for Lorentz spaces, it holds
\be
L^{2}\cdot L^{p_1,q}\;=\;L^{2,2}\cdot L^{p_1,q}\;\subset\;L^{\frac{2p_1}{2+p_1}\,,\frac{2q}{2+q}}\quad.
\ee
Since $p_1>2$, we call upon the usual Calderon-Zygmund theorem and the Sobolev-Lorentz embedding theorem of L. Tartar (Lemma 31.2 in \cite{Tar1}) to obtain that $\nabla\varphi\in W^{1,(\frac{2p_1}{2+p_1}\,,\frac{2q}{2+q})}\subset L^{p_1,\frac{2q}{2+q}}\,$ with the estimate (\ref{uuu1}) with a constant independent of $q$. \\
Next, for $1<p_0\le 2$, we use the divergence-form structure of the equation. Note that
\be\label{ww1bis}
\Delta\varphi\;=\;div\,(b\,\nabla^\perp a)\quad,
\ee
and that
\be\label{ww2bis}
b\,\nabla^\perp a\;\in\;W^{1,(p_0,q)}\cdot L^{2}\;\subset\; L^{\frac{2p_0}{2-p_0},q}\cdot L^{2,2}\;\subset\;L^{p_0,\frac{2q}{2+q}}\quad,
\ee
where we have again used the Sobolev-Lorentz embedding theorem of J.Peetre and the H\"older inequality for Lorentz spaces. The desired (\ref{uuu1}) follows immediately from (\ref{ww1bis}) and (\ref{ww2bis}).\\
%As for (\ref{u2}), it is a consequence of $p>2$ and the well-known H\"older inequality for Lorentz spaces. Namely, $L^{2,\infty}\cdot L^p\subset L^{r}$ continuously, for all $1\leq r<2p/(p+2)$. It then suffices to apply the standard Calderon-Zygmund estimate to $\Delta\varphi$ to deduce (\ref{u2}). \\
Now for $p\in (p_0,p_1)$ the estimate (\ref{uuu3}) can be reached using the Marcinkiewicz interpolation theorem just as is done in the proof of Lemma \ref{Lma1}.  \\[-3.5ex]

$\hfill\Box$ \\

\begin{Lma}\label{lm-complex-interp}
Let $a(x)$ be a positive function on the disc $D^2$ which satisfies $c_0>a(x)>c_0^{-1}$, for some positive constant $c_0$ and let $A>0$. Let $p\in[2,3]$ and define the operator
\bes
S(f)\;:=\lf[\frac{\big(A^2+a^2(x)|f|^2\big)^{(p-2)/2}}{ \big(A^2+\|f\|_p^2\big)^{(p-2)/2} }  \rg]\ f\:.
\ees
% be such that 
%\bes
%\dfrac{1}{q_2}\,<\,\dfrac{p-1}{q}\,<\,\frac{1}{q_1}\:,
%\ees
%for some $q_2> q> q_1>1$. 
If $T$ is a continuous operator from $L^s(D^2)$ to itself  for any $s\in [4/3,6]$  then, for any $\al\in (p-2,1/2)$ one has
\be
\label{A-aa}
\|TS(f)-S(Tf)\|_{p/(p-1)}\le\,2\,C\, \al^{-1}\ (p-2)\  \|T\|\ \lf[\|f\|^\al_p+\lf(A^2+\|f\|_p^2\rg)^{\al/2}\rg]\ \|f\|_p^{1-\al}\quad.
\ee
where $\|T\|:=\sup_{s\in [4/3,6]}\|T\|_{L^s\rightarrow L^s}$.
\hfill $\Box$
\end{Lma}
\noindent{\bf Proof of lemma~\ref{lm-complex-interp}.} We follow an idea from complex interpolation introduced in \cite{RW} and also used in \cite{IS}.
For $z=a+ib\in D^2_{1/2}$ we introduce
\[
S_z(f):=\lf[\frac{\big(A^2+a^2(x)|f|^2\big)^{1/2}}{ \big(A^2+\|f\|_p^2\big)^{1/2} }  \rg]^z\ f\:.
\]
Denote 
\[
p_a:=\frac{p}{1+a}\in (2\,p/3, 2p)\subset (4/3,6) \quad\mbox{ and }\quad p'_a=\frac{p_a}{p_a-1}=\frac{ p' }{1-p'a/p }\quad.
\]
It is not difficult to see that $S_z$ maps $L^p$ into $L^{p_a}$ and we shall now compute the norm of $S_z$ between these two spaces. We have
\[
\int_{D^2}|S_z(f)|^{p_a}\ dx^2\le \int_{D^2} \lf|\frac{\big(A^2+a^2(x)|f|^2\big)^{1/2}}{ \big(A^2+\|f\|_p^2\big)^{1/2} }  \rg|^{ap/(1+a)}\ |f|^{p/(1+a)}\ dx^2
\]
We first assume $\|f\|_p>A$. In this first case we have
\be
\label{A-1}
\begin{array}{l}
\ds\int_{D^2}|S_z(f)|^{p_a}\ dx^2\le C\ \|f\|_p^{-ap/(1+a)}\ \int_{D^2} \lf|\big(A^2+a^2(x)|f|^2\big)^{1/2}\rg|^{ap/(1+a)}\ |f|^{p/(1+a)}\ dx^2\\[5mm]
\ds\quad\le C\ \|f\|_p^{-ap/(1+a)}\ \int_{|f|>A} |f|^{p}\ dx^2+C\ \|f\|_p^{-ap/(1+a)}\ \int_{|f|\le A} A^{ap/(1+a)}\ |f|^{p/(1+a)}\ dx^2\quad.
\end{array}
\ee
for $a\ge 0$ first one has 
\be
\label{A-2}
\int_{|f|\le A} A^{ap/(1+a)}\ |f|^{p/(1+a)}\ dx^2\le \|f\|_p^{ap/(1+a)}\ \|f\|_p^{p/(1+a)}\quad,
\ee
and combining (\ref{A-1}) and (\ref{A-2}) one obtains
\be
\label{A-3}
\int_{D^2}|S_z(f)|^{p_a}\ dx^2\le C\ \|f\|_p^{p_a}\quad.
\ee
Nor for $a<0$, still under the assumption that $\|f\|_p>A$ one has
\be
\label{A-4}
\int_{|f|\le A} A^{ap/(1+a)}\ |f|^{p/(1+a)}\ dx^2\le \int_{|f|\le A}|f|^p\ dx^2\quad.
\ee
Combining (\ref{A-1}) and (\ref{A-4}) gives again (\ref{A-3}). We assume now $\|f\|_p\le A$. In this case we have
\be
\label{A-5}
\begin{array}{l}
\ds\int_{D^2}|S_z(f)|^{p_a}\ dx^2\le C\ A^{-ap/(1+a)}\ \int_{D^2} \lf|\big(A^2+a^2(x)|f|^2\big)^{1/2}\rg|^{ap/(1+a)}\ |f|^{p/(1+a)}\ dx^2\\[5mm]
\ds\quad\le C\ A^{-ap/(1+a)}\ \int_{|f|>A} |f|^{p}\ dx^2+C\ A^{-ap/(1+a)}\ \int_{|f|\le A} A^{ap/(1+a)}\ |f|^{p/(1+a)}\ dx^2\quad.
\end{array}
\ee
First, for $a>0$, under this assumption  $\|f\|_p\le A$ one gets
\be
\label{A-6}
A^{-ap/(1+a)}\ \int_{|f|>A} |f|^{p}\ dx^2\le \|f\|_{p}^{-ap/(1+a)}\ \|f\|_p^p=\|f\|_p^{p_a}\quad.
\ee
and in this case again we have (\ref{A-3}). Now, for $a<0$ and $\|f\|_p\le A$ we have from (\ref{A-5})
\be
\label{A-7}
\int_{D^2}|S_z(f)|^{p_a}\ dx^2\le C\ \lf(\frac{A}{\|f\|_p}\rg)^{-ap/(1+a)}\ \|f\|_p^{p_a}\quad.
\ee
This gives in all cases
\be
\label{A-8}
\|S_z(f)\|_{p_a}\le C\ \lf[1+\lf(\frac{\|f\|^2_p}{A^2+\|f\|_p^2}\rg)^{a/2}\rg]\ \|f\|_p\quad.
\ee
Denote
\[
Q_z(g):= \lf(\frac{\|g\|_{p'}}{|g|}\rg)^{\frac{\ov{z}p'}{p}}\ g\quad.
\]
Observe that we have
\be
\label{A-9}
\|Q_z(g)\|_{p'_a}= \|g\|_{p'}
\ee
Let $f\in L^p(D^2)$ and $g\in L^{p'}(D^2)$ and let $T$ be a continuous operator from $L^s(D^2)$ into $L^s(D^2)$ for any $s\in[4/3,6]$ such that
\[
\|T\|:=\sup_{s\in [4/3,6]}\||T\||_{L^s\rightarrow L^s}<+\infty\quad.
\]
 We introduce the function defined on the disc $D^2$ by
\[
\varphi(z):=\int_{D^2}\lf[TS_z(f)-S_z(Tf)\rg] \ov{Q_z(g)}\ dx^2\quad.
\]
$\varphi$ is obviously holomorphic and one has
\be
\label{A-10}
|\varphi(z)|\le\ \|T\|\ \|S_z(f)\|_{p_a}\ \|Q_z(g)\|_{p'_a}\le  C\  \|T\|\ \lf[1+\lf(\frac{\|f\|^2_p}{A^2+\|f\|_p^2}\rg)^{\Re(z)/2}\rg]\ \|f\|_p\ \|g\|_{p'}\quad.
\ee
Hence
\be
\label{A-11}
\sup_{|z|\le \al}|\varphi(z)|\le  C\  \|T\|\ \lf[1+\lf(\frac{A^2+\|f\|_p^2}{\|f\|^2_p}\rg)^{\al/2}\rg]\ \|f\|_p\ \|g\|_{p'}
\ee
Using Schwartz lemma we obtain the pointwise estimate on $D^2_{\al}$
\be
\label{A-12}
|\varphi(z)|\le\ C\ \al^{-1}\ |z|\  \|T\|\ \lf[1+\lf(\frac{A^2+\|f\|_p^2}{\|f\|^2_p}\rg)^{\al/2}\rg]\ \|f\|_p\ \|g\|_{p'}
\ee
Hence, using (\ref{A-9}), after observing that for $g=(h\, \|h\|_{p'_a}^{-\ov{z}\,p'/p})^{1/(1-\ov{z}p'/p)}$ one has $Q_z(g)=h$, we deduce that
\be
\label{A-13}
\sup_{\|h\|_{p'_a}\le 1}\int_{D^2}\lf[TS_z(f)-S_z(Tf)\rg] \ov{h}\ dx^2\le\,2\,C\, \al^{-1}\ |z|\  \|T\|\ \lf[1+\lf(\frac{A^2+\|f\|_p^2}{\|f\|^2_p}\rg)^{\al/2}\rg]\ \|f\|_p\ 
\ee
We apply (\ref{A-13}) to $z=p-2$ and we obtain
\be
\label{A-14}
\|TS_{p-2}(f)-S_{p-2}(Tf)\|{p/(p-1)}\le\,2\,C\, \al^{-1}\ |z|\  \|T\|\ \lf[\|f\|^\al_p+\lf(A^2+\|f\|_p^2\rg)^{\al/2}\rg]\ \|f\|_p^{1-\al}\quad.
\ee
which implies the lemma.\hfill $\Box$
\begin{Lma}
\label{lma-A-6}
Let $p\in [2,3)$, $r<1$, $B\in {\R}$, $a(x)\in L^\infty(D^2_r)$ and $\vec{n}$ such that $\nabla\vec{n}\in L^p(D^2_r)$ and assume  $c_0^{-1}\le a(x)\le c_0$. Denote $\zeta\in L^{p/(p-1)}$
to be the solution of
\be
\label{A-15}
\lf\{
\begin{array}{l}
\ds-\Delta\vec{\zeta}= \mbox{div}\lf(  \lf(B^2+a^2(x)\, |\nabla^\perp\vec{n}|^2\rg)^{p/2-1}\ \nabla^\perp\vec{n}  \rg)   \quad\quad\mbox{ in } D_r^2 \\[5mm]
\ds \vec{\zeta}=0\quad\quad\mbox{ on }\p D_r^2
\end{array}
\rg.
\ee
Then the following inequality holds for any $p-2<\al<1/2$
\be
\label{A-16}
\|\nabla\vec{\zeta}\|_{L^{p/(p-1)}(D^2_r)}\le  C\  \frac{(p-2)}{\al}\ \lf(\|\nabla\vec{n}\|^{p-1}_{L^p(D^2_r)}+(r^{2/p}\,B)^{(p-2)+\al}\, \|\nabla\vec{n}\|^{1-\al}_{L^p(D^2_r)}\rg)\quad,
\ee
where $C$ is independent of $p$, $\al$, $r$, $B$ and $\vec{n}$.
\hfill $\Box$
\end{Lma}
\noindent{\bf Proof of lemma~\ref{lma-A-6}.} Let $\vec{n}_r(x):=\vec{n}(r\,x)$ and $\vec{\zeta}_r(x):=\vec{\zeta}(r\,x)$. The following system is satisfied
\be
\label{A-17}
\lf\{
\begin{array}{l}
\ds-\,\Delta\vec{\zeta}_r= \,r^{-p+2}\mbox{div}\lf(  \lf(B^2\, r^2+a^2(x)\, |\nabla^\perp\vec{n}_r|^2\rg)^{p/2-1}\ \nabla^\perp\vec{n}_r  \rg)   \quad\quad\mbox{ in } D_1^2 \\[5mm]
\ds \vec{\zeta}_r=0\quad\quad\mbox{ on }\p D_1^2
\end{array}
\rg.
\ee
Let $T$ be the operator which to an $L^s(D^2)$ vector Field $X$ assigns $\nabla u$ the solution to
\[
\begin{array}{l}
-\mbox{div}(\nabla u)= \mbox{div}(X) \quad\quad\mbox{ on }\p D^2\\[5mm]
u=0\quad\quad\mbox{ on }\p D^2
\end{array}
\]
Calderon Zygmund Theory gives that $T$ is continuous from $L^s(D^2)$ into $L^s(D^2)$ and $\sup_{s\in [4/3,6]}\|T\|_{L^s\rightarrow L^s}=\|T\|<+\infty$. Observe that
\be
\label{A-18}
T\lf(\lf(B^2\, r^2+a^2(x)\, |\nabla^\perp\vec{n}_r|^2\rg)^{p/2-1}\ \nabla^\perp\vec{n}_r \rg)=r^{p-2}\, \nabla\,\vec{\zeta}_r\quad\mbox{ and }\quad T(\nabla^\perp\vec{n}_r)=0
\ee
Denote $A:=B\, r$ and 
\[
S(f):=\lf(\frac{A^2+a^2(x)\, |f|^2}{A^2+\|f\|^2_p}\rg)^{p/2-1}\ f
\]
Combining (\ref{A-aa}) with (\ref{A-18}) one obtains
\be
\label{A-19}
\frac{r^{p-2}\, \|\nabla\vec{\zeta}_r\|_{L^{p/(p-1)}(D^2_1)}}{\lf(A^2+\|\nabla\vec{n}_r\|^2_p\rg)^{p/2-1}}\le  2\,C\, \al^{-1}\ (p-2)\  \|T\|\ \lf[\|\nabla\vec{n}_r\|^\al_p+\lf(A^2+\|\nabla\vec{n}_r\|_p^2\rg)^{\al/2}\rg]\ \|\nabla\vec{n}_r\|_p^{1-\al}\ \quad.
\ee
Observe that 
\[
\|\nabla\vec{n}_r\|_{L^p(D^2_1)}=r^{1-2/p}\ \|\nabla\vec{n}\|_{L^p(D^2_r)}\quad\mbox{ and }\quad \|\nabla\vec{\zeta}_r\|_{L^{p/(p-1)}(D^2_1)}=r^{-1+2/p}\ \|\nabla\vec{\zeta}\|_{L^{p/(p-1)}(D^2_r)}\quad.
\]
Then
\be
\label{A-20}
\,\|\nabla\vec{\zeta}\|_{L^{p/(p-1)}(D^2_r)}\le  2\,C\  \frac{(p-2)}{\al}\  \|T\|\ \lf(\|\nabla\vec{n}\|^{p-1}_{L^p(D^2_r)}+(r^{2/p}\,B)^{(p-2)+\al}\, \|\nabla\vec{n}\|^{1-\al}_{L^p(D^2_r)}\rg)
\ee
which implies the lemma.\hfill $\Box$

\end{document}